\input german
\newwrite\num
\catcode`\^^Z=9
\catcode`\^^M=10
\input german
\output={\if N\header\headline={\hfill}\fi
\plainoutput\global\let\header=Y}
\magnification\magstep1
\tolerance = 500
\hsize=14.4true cm
\vsize=22.5true cm
\parindent=6true mm\overfullrule=2pt
\newcount\kapnum \kapnum=0
\newcount\parnum \parnum=0
\newcount\procnum \procnum=0
\newcount\nicknum \nicknum=1
\font\ninett=cmtt9

\font\ninebf=cmbx9

\font\sixbf=cmbx6
\font\ninesl=cmsl9

\font\nineit=cmti9

\font\ninerm=cmr9

\font\sixrm=cmr6
\font\ninei=cmmi9
\font\eighti=cmmi8
\font\sixi=cmmi6
\skewchar\ninei='177 \skewchar\eighti='177 \skewchar\sixi='177
\font\ninesy=cmsy9
\font\eightsy=cmsy8
\font\sixsy=cmsy6
\skewchar\ninesy='60 \skewchar\eightsy='60 \skewchar\sixsy='60
\font\titelfont=cmr10 scaled 1440
\font\paragratit=cmbx10 scaled 1200

\font\name=cmcsc10
\font\emph=cmbxti10

\font\tenmsbm=msbm10
\font\sevenmsbm=msbm7
%

%
\font\got=eufm10

\font\teneufm=eufm10
\font\seveneufm=eufm7
\font\fiveeufm=eufm5
\newfam\eufmfam
\textfont\eufmfam=\teneufm
\scriptfont\eufmfam=\seveneufm
\scriptscriptfont\eufmfam=\fiveeufm

\font\tenmsam=msam10
\font\sevenmsam=msam7
\font\fivemsam=msam5
\newfam\msamfam
\textfont\msamfam=\tenmsam
\scriptfont\msamfam=\sevenmsam
\scriptscriptfont\msamfam=\fivemsam
\font\tenmsbm=msbm10
\font\sevenmsbm=msbm7
\font\fivemsbm=msbm5
\newfam\msbmfam
\textfont\msbmfam=\tenmsbm
\scriptfont\msbmfam=\sevenmsbm
\scriptscriptfont\msbmfam=\fivemsbm
\def\Bbb#1{{\fam\msbmfam\relax#1}}
\def\cz{{\kern0.4pt\Bbb C\kern0.7pt}
}
\def\ez{{\kern0.4pt\Bbb E\kern0.7pt}
}
\def\fz{{\kern0.4pt\Bbb F\kern0.3pt}}
\def\gz{{\kern0.4pt\Bbb Z\kern0.7pt}}
\def\hz{{\kern0.4pt\Bbb H\kern0.7pt}
}
\def\kz{{\kern0.4pt\Bbb K\kern0.7pt}
}
\def\nz{{\kern0.4pt\Bbb N\kern0.7pt}
}
\def\oz{{\kern0.4pt\Bbb O\kern0.7pt}
}
\def\rz{{\kern0.4pt\Bbb R\kern0.7pt}
}
\def\sz{{\kern0.4pt\Bbb S\kern0.7pt}
}
\def\pz{{\kern0.4pt\Bbb P\kern0.7pt}
}
\def\qz{{\kern0.4pt\Bbb Q\kern0.7pt}
}
\newskip\ttglue
\def\ninepoint{\def\rm{\fam0\ninerm}%
  \textfont0=\ninerm \scriptfont0=\sixrm \scriptscriptfont0=\fiverm
  \textfont1=\ninei \scriptfont1=\sixi \scriptscriptfont1=\fivei
  \textfont2=\ninesy \scriptfont2=\sixsy \scriptscriptfont2=\fivesy
  \textfont3=\tenex \scriptfont3=\tenex \scriptscriptfont3=\tenex
  \def\it{\fam\itfam\nineit}%
  \textfont\itfam=\nineit
  \def\sl{\fam\slfam\ninesl}%
  \textfont\slfam=\ninesl
  \def\bf{\fam\bffam\ninebf}%
  \textfont\bffam=\ninebf \scriptfont\bffam=\sixbf
   \scriptscriptfont\bffam=\fivebf
  \def\tt{\fam\ttfam\ninett}%
  \textfont\ttfam=\ninett
  \tt \ttglue=.5em plus.25em minus.15em
  \normalbaselineskip=11pt
  \font\name=cmcsc9
  \let\sc=\sevenrm
  \let\big=\ninebig
  \setbox\strutbox=\hbox{\vrule height8pt depth3pt width0pt}%
  \normalbaselines\rm
  \def\sl{\it}}

\headline={\ifodd\pageno\rightheadline\else\leftheadline\fi}
\def\rightheadline{\ninepoint Paragraphen"uberschrift\hfill\folio}
\def\leftheadline{\ninepoint\folio\hfill Chapter"uberschrift}
\let\header=Y
\def\titel#1{\need 9cm \vskip 2truecm
\parnum=0\global\advance \kapnum by 1
{\baselineskip=16pt\lineskip=16pt\rightskip0pt
plus4em\spaceskip.3333em\xspaceskip.5em\pretolerance=10000\noindent
\titelfont Chapter \uppercase\expandafter{\romannumeral\kapnum}.
#1\vskip2true cm}\def\leftheadline{\ninepoint
\folio\hfill Chapter \uppercase\expandafter{\romannumeral\kapnum}.
#1}\let\header=N
}
\def\Titel#1{\need 9cm \vskip 2truecm
\global\advance \kapnum by 1
{\baselineskip=16pt\lineskip=16pt\rightskip0pt
plus4em\spaceskip.3333em\xspaceskip.5em\pretolerance=10000\noindent
\titelfont\uppercase\expandafter{\romannumeral\kapnum}.
#1\vskip2true cm}\def\leftheadline{\ninepoint
\folio\hfill\uppercase\expandafter{\romannumeral\kapnum}.
#1}\let\header=N
}
\def\need#1cm {\par\dimen0=\pagetotal\ifdim\dimen0<\vsize
\global\advance\dimen0by#1 true cm
\ifdim\dimen0>\vsize\vfil\eject\noindent\fi\fi}
\def\neupara#1{\par\penalty-2000
\procnum=0\global\advance\parnum by 1
\vskip1cm\noindent{\paragratit \the\parnum. #1}%
\def\rightheadline{\ninepoint\S\the\parnum.\ #1\hfill \folio}%
\vskip 8mm\noindent}
\def\Proclaim #1 #2\finishproclaim {\bigbreak\noindent
{\bf#1\unskip{}. }{\it#2}\medbreak\noindent}
%
\gdef\proclaim #1 #2 #3\finishproclaim {\bigbreak\noindent%
\global\advance\procnum by 1
{%
{\relax\ifodd \nicknum
\hbox to 0pt{\vrule depth 0pt height0pt width\hsize
   \quad \ninett#3\hss}\else {}\fi}%
\bf\the\parnum.\the\procnum\ #1\unskip{}. }
{\it#2}
\immediate\write\num{\string\def
 \expandafter\string\csname#3\endcsname
 {\the\parnum.\the\procnum}}
\medbreak\noindent}
\newcount\stunde \newcount\minute \newcount\hilfsvar
\def\uhrzeit{
    \stunde=\the\time \divide \stunde by 60
    \minute=\the\time
    \hilfsvar=\stunde \multiply \hilfsvar by 60
    \advance \minute by -\hilfsvar
    \ifnum\the\stunde<10
    \ifnum\the\minute<10
    0\the\stunde:0\the\minute~Uhr
    \else
    0\the\stunde:\the\minute~Uhr
    \fi
    \else
    \ifnum\the\minute<10
    \the\stunde:0\the\minute~Uhr
    \else
    \the\stunde:\the\minute~Uhr
    \fi
    \fi
    }

 \def\bfH{{\bf H}}

\def\calA{{\cal A}} \def\calB{{\cal B}}

 \def\calH{{\cal H}}
\def\calI{{\cal I}} 
 
\def\calM{{\cal M}} 
 
 \def\calR{{\cal R}}

\def\gotx{\hbox{\got x}} 
\def\goty{\hbox{\got y}}

\def\Aut{\mathop{\rm Aut}\nolimits}

\def\dim{\mathop{\rm dim}\nolimits}

\def\GL{\mathop{\rm GL}\nolimits}

\def\Hom{\mathop{\rm Hom}\nolimits}

\def\im{\mathop{\rm Im}\nolimits} \def\Im{\im}

\def\kernel{\mathop{\rm kernel}\nolimits}

\def\mod{\mathop{\rm mod}\nolimits}
\def\O{{\rm O}}
\def\U{{\rm U}}

\def\proj{\mathop{\rm proj}\nolimits}

\def\re{\mathop{\rm Re}\nolimits}
\def\Re{\re}

\def\SL{\mathop{\rm SL}\nolimits}

\def\Sp{\mathop{\rm Sp}\nolimits}

\def\Sym{\mathop{\rm Sym}\nolimits}
\def\boxit#1{
  \vbox{\hrule\hbox{\vrule\kern6pt
  \vbox{\kern8pt#1\kern8pt}\kern6pt\vrule}\hrule}}
\def\Boxit#1{
  \vbox{\hrule\hbox{\vrule\kern2pt
  \vbox{\kern2pt#1\kern2pt}\kern2pt\vrule}\hrule}}

\def\zwischen#1{\bigbreak\noindent{\bf#1\medbreak\noindent}}

\def\smallni{\smallskip\noindent }
\def\medni{\medskip\noindent }
\def\bigni{\bigskip\noindent }
\def\Isom{\mathop{\;{\buildrel \sim\over\longrightarrow }\;}}
\def\lo{\longrightarrow}

\def\loma{\longmapsto}
\def\pii{\pi {\rm i}}

\def\set#1{\bigl\{\,#1\,\bigr\}}

\def\square{\hbox{\hbox to 0pt{$\sqcup$\hss}\hbox{$\sqcap$}}}
\def\qed{\ifmmode\square\else{\unskip\nobreak\hfil
\penalty50\hskip3em\null\nobreak\hfil\square
\parfillskip=0pt\finalhyphendemerits=0\endgraf}\fi}
\def\pn{\the\parnum.\the\procnum}
\def\downmapsto{{\buildrel
        {\vbox{\hbox{\hskip.2pt$\scriptstyle-$}}}
        \over{\raise7pt\vbox{\vskip-4pt\hbox{$\textstyle\downarrow$}}}}}

\let\header=N
\newwrite\num
\openout\num=eacht6.num

\def\Mp{{\rm Mp}}

\def\imag{{{\rm i}_1}}

\def\trace{\hbox{\rm trace}}

\def\Lai{2.1}
\def\RbC{2.2}
\def\LRl{2.3}

\def\Dsi{2.5}
\def\Kcs{2.6}
\def\SsR{2.7}
\def\Qrd{2.8}
\def\TbK{2.9}
\def\QcQ{2.10}

\def\KF{3.2}
\def\Pfc{3.3}

\def\Nvv{4.7}
\def\As{5.1}
\def\PIo{5.2}
\def\Sev{5.3}
\def\Lmm{6.1}

\def\Pfs{6.3}
\def\Sta{6.4}

\def\Tab{6.6}
\def\KR{7.1}

\def\PB{7.3}

\def\BaB{8.2}
\def\KU{8.3}

\def\NIk{9.1}

\def\SAz{9.3}
\def\Smz{9.4}

\def\Kgn{10.1}
\def\nSa{10.2}

\def\Fst{10.4}
\def\cDA{10.5}
\def\LKa{10.6}
\def\LKb{10.7}
\def\Sdc{10.8}
\def\RFb{10.9}

\def\RAND#1{\vskip0pt\hbox to 0mm{\hss\vtop to 0pt{%
  \raggedright\ninepoint\parindent=0pt%
  \baselineskip=1pt\hsize=2cm #1\vss}}\noindent}
\noindent\centerline{\titelfont Modular forms for the}
\centerline{\titelfont even unimodular
lattice of signature  (2,10)}%
\vskip 1.5cm
\leftline{\it \hbox to 6cm{Eberhard Freitag\hss}
Riccardo Salvati
Manni  }
  \leftline {\it  \hbox to 6cm{Mathematisches Institut\hss}
Dipartimento di Matematica, }
\leftline {\it  \hbox to 6cm{Im Neuenheimer Feld 288\hss}
Piazzale Aldo Moro, 2}
\leftline {\it  \hbox to 6cm{D69120 Heidelberg\hss}
 I-00185 Roma, Italy. }
\leftline {\tt \hbox to 6cm{freitag@mathi.uni-heidelberg.de\hss}
salvati@mat.uniroma1.it}
\def\text{\hbox}

\def\trace{{\rm trace}}
\def\leftheadline{\ninepoint\folio\hfill
Modular forms}%
\def\rightheadline{\ninepoint Introduction\hfill \folio}%
\headline={\ifodd\pageno\rightheadline\else\leftheadline\fi}
\bigni
\centerline{\paragratit \rm  2005}%
\vskip5mm\noindent%
\let\header=N%
\def\imag{{\rm i}}%
\def\Mp{{\rm Mp}}%
{\paragratit Introduction}%
\medskip\noindent
More than forty years ago, Igusa proved in [Ig] a ``Fundamental Lemma'',
according to which, theta constants give an explicit  generically
injective map from some modular varieties related to the symplectic
group $\Sp(n,\rz)$ into the  projective space. The embedded varieties
satisfy some quartic relations, the  so called Riemann's relations.
We expect that similar results hold for orthogonal instead of symplectic groups,
where one should use typical ``orthogonal constructions'' instead of theta series.
\smallskip
Some years ago, Borcherds described in [Bo1] two  methods for
constructing modular forms on
modular varieties related to the orthogonal group $\O (2,n)$. They
are the so called Borcherds' additive  and multiplicative lifting.
The multiplicative lifting has been used by Borcherds himself and
other authors to construct modular forms  with known vanishing locus
and interesting properties.
The additive lifting has been used to
construct explicit maps from
some modular varieties
related to $\O(2,4)$, $\O(2,6)$, $\O(2,8)$, $\O(2,10)$
and also  for some unitary groups as $\U(1,4)$ and $\U(1,5)$.
cf. [FH], [AF], [Fr2], [Fr3], [FS], [Ko1], [Ko2].
\smallskip
In this paper we try  a more systematic treatment in certain level 2 cases:
We are mainly interested in even lattices $L$ of signature $(2,n)$, such that
the discriminant $L'/L$ is a vector space over the field of two elements and such
the the induced quadratic form $q:L'/L\to\fz_2$ is of even type.
This means that the dimension $2m$ of $L'/L$ is even and that $q$ is of the
form
$$q(x)=x_1x_2+\cdots+x_{2m-1}x_{2m}$$
with respect to a suitable basis.
This implies also that $n\equiv 2$ mod 8.
To every such finite geometry $(\fz_2^{2m},q)$
we attach a certain graded algebra $\calR_m/\calI_m$.
We expect that this algebra is related to the structure of certain
modular varieties. There are three cases where this can be verified.
The first one has been investigated by
Kondo's in connection with his research about the
moduli space of Enriques surfaces. The lattice here is
$U\oplus\sqrt2U\oplus(-\sqrt 2 E_8)$ and in this case $m=5$.
In section 11 we describe in some detail this case again,
giving the connection with our results and correcting
some error in Kondo's paper.
\smallskip
Another case, related to the configuration space of 8 points
in the projective line, appears in [Ko2] and [Koi]. It is related to the lattice
$U\oplus\sqrt2U\oplus(-D_4)\oplus(-D_4)$, i.e. we have $m=3$ in this case.
\smallskip
A third member is studied in this paper, namely the lattice
$\sqrt 2(U\oplus U\oplus(-E_8))$. Here $m=6$. In all three cases the
ideal $\calI_m\subset\calR_m$ occurs. The modular variety, which we investigate,
is a covering of Kondo's variety
(modular variety of Enriques surfaces), which can be treated in a similar manner.
But we have a bigger symmetry group $\O(\fz_2^{12})$ instead of
$\O(\fz_2^{10})$ with the effect that some difficulties of Kondo's approach
disappear.
\smallskip
The ideal $\calI_m$ is defined by quadratic relations in a polynomial ring
 $\calR_m$. Kondo considers instead of this only quartic relations.
In [Ko2] and [Koi] has been explained
in the case $m=3$ that the quartic relations are consequences
of quadratic relations. This is true for arbitrary $m$, hence also in Kondo's
$\O(2,10)$-case.
\smallskip
\neupara{Lattices}%
A lattice $L$ is a free abelian group together with a
symmetric non degenerated bilinear quadratic form
$$q_L:L\lo\qz.$$
The associated bilinear form is $(x,y)=q_L(x,y)-q_L(x)-q_L(y)$.
We extend $(x,y)$ to a $\rz$-bilinear from on $V=L\otimes_\gz\rz$
and moreover to a $\cz$-bilinear form on $L\otimes_\gz\cz$
and we use the notation $q_L(x)=(x,x)/2$ also in this situation.
\smallskip
The  {\it norm\/} of a vector $a\in V$ is $(a,a)=2q_L(a)$.
The lattice is called even if $q_L(a)$ is integral for all $a$.
The dual lattice $L'\subset L\otimes_\gz\rz$ consists of all
elements $a$ such that $(a,x)$ is integral for all $x\in L$.
The lattice is called unimodular if $L=L'$. We recall that
an even unimodular lattice of signature $(m,n)$ exists if and only if
$m-n$ is divisible by 8 and that it is determined up to isomorphism
in this case. This lattice is denoted by II$_{m,n}$.
We use the usual notation $U=\hbox{II}_{1,1}$.
\smallskip
In this paper we will treat mainly   even lattices, such that $2L'\subset L$.
Let us assume that $V$ is of signature $(2,n)$.
We recall that there is a unique subgroup of index two
$\O^+(V)\subset\O(V)$, which contains all reflections
$$x\loma x-{(a,x)\over q_L(a)}a$$
for vectors $a\in V$ with {\it negative\/} $q_L(a)$.
We also consider the integral orthogonal group
$\O(L)$ and
$$\O^+(L)=\O(L)\cap\O^+(V).$$
Since $L$ is even, $q_L:L'\to\qz$ factors through a map
$$\bar q_L:L'/L\lo\qz/\gz,$$
which is called the discriminant of $(L,q_L)$.
\smallskip
We denote by
$\Aut(A)$ the group of automorphisms of an abelian group and by
$$\O(L'/L)=\O(L'/L,\bar q_L)\subset\Aut(L'/L)$$
the subgroup which preserves $\bar q_L$.
The kernel of the map
$\O(L)\lo\O(L'/L)$
is called the discriminant kernel.
We  define
$$\Gamma_L=\hbox{kernel}\bigl(\O^+(L)\lo\O(L'/L)\bigr).$$
Let $\cal H$ be the bounded symmetric domain  of type IV associated to $\O(V)$,
 then we shall consider the modular variety
 $\calH/\Gamma_L$.

\neupara{Finite geometries and related projective varieties}%
We consider $\fz_2^{2m}$ equipped with the quadratic form
$$q(x)=x_1x_2+\cdots+x_{2m-1}x_{2m}$$
and the associated bilinear form
$$(x,y)=q(x+y)-q(x)-q(y).$$
The orthogonal group (group of all automorphims which preserve $q$)
is denoted by $\O(\fz_2^{2m})$.
Its order is
$$2^{m(m-1)+1}(2^m-1)\prod_{i=1}^{m-1}(2^{2i}-1).$$
The number of isotropic elements is $2^{m-1}(2^m+1)$ and the number of
anisotropic elements is $2^{m-1}(2^m-1)$.
We will need the following refined information:
\proclaim
{Lemma}
{1) Let $a\in\fz_2^{2m}$ be a non-zero isotropic element.
\smallskip
The number
of all non-zero isotropic elements $b$ with $(a,b)=0$ is
$$2^{m-1}(2^{m-1}+1)-1$$
and the number
of all anisotropic elements $b$ with $(a,b)= 0$ is
$$2^{m-1}(2^{m-1}-1).$$
The number
of all isotropic elements with $(a,b)\ne 0$ is
$$2^{2(m-1)}.$$
and the number
anisotropic elements $b$ with $(a,b)\ne 0$ is also
$$2^{2(m-1)}.$$
2) Let $a\in\fz_2^{2m}$ be an anisotropic element.
\smallskip
The number
of all non-zero isotropic elements $b$ with $(a,b)=0$ is
$$2^{2(m-1)}-1$$
and the number
of all anisotropic elements $b$ with $(a,b)= 0$ is
$$2^{2(m-1)}.$$
The number
of all isotropic elements with $(a,b)\ne 0$ is
$$2^{m-1}(2^{m-1}+1).$$
and the number
anisotropic elements $b$ with $(a,b)\ne 0$ is also
$$2^{m-1}(2^{m-1}-1).$$
}
Lai%
\finishproclaim
A sub-space $A\subset\fz_2^{2m}$ is called totally isotropic if all elements
are isotropic.
The maximal totally isotropic spaces are of dimension $m$. They form one orbit under
the orthogonal group.
\proclaim
{Lemma}
{Let $A\subset\fz_2^{2m}$ be a totally isotropic subspace of dimension $m-2$.
There are
precisely 6 maximal totally isotropic spaces (of dimension $m$)
which contain $A$.
}
RbC%
\finishproclaim
{\it Proof.\/}
We decompose $\fz_2^{2m}=\fz_2^{2m-4}\oplus\fz_2^4$. Let $A\subset\fz_2^{2m-4}$
be a maximal totally isotropic subspace. Then $A$ can be consider as a
$(m-2)$-dimensional totally isotropic space in $\fz_2^{2m}$. The maximal
totally isotropic subspaces which contain  it are of the form $A\oplus B$ where
$B\subset\fz_2^4$ is a two dimensional totally isotropic subspace of $\fz_2^4$.
It is easy to check that there are 6 such spaces.\qed
\smallskip
We list the non-zero elements of these subspaces $B$:
\vskip2mm
\halign{\quad#\hfil&\quad#\hfil&\quad#\hfil&\quad#\hfil\cr
$I_1$&(1,0,0,0)&(0,0,1,0)&(1,0,1,0)\cr
$I_2$&(1,0,0,0)&(0,0,0,1)&(1,0,0,1)\cr
$I_3$&(0,1,0,0)&(0,0,0,1)&(0,1,0,1)\cr
$I_4$&(0,1,0,0)&(0,0,1,0)&(0,1,1,0)\cr
$I_5$&(1,0,0,1)&(0,1,1,0)&(1,1,1,1)\cr
$I_6$&(1,0,1,0)&(0,1,0,1)&(1,1,1,1)\cr
}
\smallskip\noindent
This table shows:
\proclaim
{Lemma}
{Let $A\subset\fz_2^{2m}$ be a totally isotropic subspace of dimension $m-2$.
The space which is generated by the characteristic functions
of the maximal (=m-dimensional) totally isotropic
subspaces containing $A$ has dimension $5$.
If $I_1,\dots,I_6$ denotes these subspaces in a suitable ordering
and $\chi_1,\dots,\chi_6$ are their characteristic functions
then one has the relation
$$\chi_{1}+\chi_{3}+\chi_{5}=\chi_{2}+\chi_{4}+\chi_{6}$$
and furthermore the nine spaces
$I_i\cap I_j$ with odd $i$ and even $j$ are\/
$(m-1)$-dimensional. }
LRl%
\finishproclaim
{\it Proof.\/} This is a calculation inside $\fz_2^4$ which can be done
by means of the above table.\qed\smallskip

Recall that the dimension of a maximal totally isotropic subspace of $\fz_2^{2m}$ is
$m$. A typical example is the subspace defined by $x_2=x_4=\cdots= x_{2m}=0$.
We consider now $(m-1)$-dimensional totally isotropic subspaces $N\subset\fz_2^{2m}$,
for example that one which is defined by the additional condition $x_{2m-1}=0$.
\proclaim
{Lemma}
{Let $N\subset\fz_2^{2m}$ be  a $(m-1)$-dimensional totally isotropic subspace.
There exist
exactly two cosets different from $N$, which contain only isotropic vectors
and one coset which contains only anistropic vectors.}
Oai%
\finishproclaim
{\it Proof.\/} Because we have a transitive action of the orthogonal group
on the set of $(m-1)$-dimensional totally isotropic subspaces, it is sufficient to treat
the above example. In this case the  relevant orbits are represented vectors
of the form $(1,0,1,0,1,0,1,0,1,0,*,*)$.\qed
\proclaim
{Definition}
{A {\emph star} is a set of $2^{m-1}$ anisotropic vectors\/ of $\fz_2^{2m}$
which form a coset of a $(m-1)$- dimensional totally isotropic space.}
Dsi%
\finishproclaim
This isotropic space is determined by the star.
A simple calculation shows:
\proclaim {Lemma} {Let $M\subset\fz_2^{2m}$ be
star with the underlying $(m-1)$-dimensional totally isotropic subspace
$N\subset\fz_2^{2m}$. There exist precisely two m-dimensional
(i.e.\ maximal) totally isotropic subspaces  which contain $N$. The
difference of their characteristic functions is an element of
$\cz[\fz_2^{2m}]$ with the following property. It
changes the sign if one applies the reflections $x\mapsto x+(a,x)a$ along the\/
$2^{m-1}$ anisotropic vectors $a$ of the associated star.}
Kcs%
\finishproclaim
Next we consider all $(m-1)$-dimensional totally isotropic
subspaces which contain our $(m-2)$-dimensional standard space
$A$. They are of the form $A+\fz_2\alpha$, where $\alpha$ is a
non-zero isotropic. Hence there are nine such spaces. We list them
and the corresponding star (=unique orbit of this space which
contains only anisotropic elements):
We only list the $\fz_2^4$-%
components. The first $2m-4$ components vary arbitrarily in $A$.
\vskip 5mm
\vbox{\halign{\quad#\hfil&\quad\qquad#\hfil&\qquad#\hfil&\quad #\hfil&\quad #\hfil\cr
\hbox to 0pt{totally isotropic spaces\hss}
&&\hbox to 0pt{corresponding star\hss}\cr
\noalign{\smallskip}
$I_1\cap I_2$&(1,0,0,0)&$O_1$&(0,0,1,1)&(1,0,1,1)\cr
$I_3\cap I_4$&(0,1,0,0)&$O_2$&(0,0,1,1)&(0,1,1,1)\cr
$I_1\cap I_4$&(0,0,1,0)&$O_3$&(1,1,0,0)&(1,1,1,0)\cr
$I_3\cap I_2$&(0,0,0,1)&$O_4$&(1,1,0,0)&(1,1,0,1)\cr
$I_1\cap I_6$&(1,0,1,0)&$O_5$&(0,1,1,1)&(1,1,0,1)\cr
$I_5\cap I_2$&(1,0,0,1)&$O_6$&(0,1,1,1)&(1,1,1,0)\cr
$I_5\cap I_4$&(0,1,1,0)&$O_7$&(1,0,1,1)&(1,1,0,1)\cr
$I_3\cap I_6$&(0,1,0,1)&$O_8$&(1,0,1,1)&(1,1,1,0)\cr
$I_5\cap I_6$&(1,1,1,1)&$O_9$&(1,1,0,0)&(0,0,1,1)\cr
}}
\medni
Recall that every $(m-1)$-dimensional totally isotropic space $B$ is associated to a
star, which is the unique coset consisting of anistropic elements.
There are two $m$-dimensional totally isotropic spaces which contain $B$.
\proclaim{Remark}
{Let $\psi_i$ be the characteristic functions of the nine stars above.
They satisfy the relations
$$\eqalign{
\psi_1+\psi_3=\psi_8+\psi_9\cr
\psi_1+\psi_4=\psi_7+\psi_9\cr
\psi_1+\psi_5=\psi_2+\psi_7\cr}
\quad\eqalign{
\psi_1+\psi_6=\psi_2+\psi_8\cr
\psi_2+\psi_3=\psi_6+\psi_9\cr
\psi_3+\psi_4=\psi_5+\psi_9\cr}
\quad\eqalign{
\psi_3+\psi_5=\psi_4+\psi_6\cr
\psi_3+\psi_7=\psi_4+\psi_8\cr
\psi_5+\psi_8=\psi_6+\psi_7\cr}

$$}
SsR%
\finishproclaim
This relations are related to certain quadratic
relations of characteristic functions.
\proclaim
{Lemma}
{Let $A\subset\fz_2^{2m}$ be an totally isotropic subspace of
dimension $m-2$ and
$I_1,\dots,I_6$ the maximal totally isotropic subspaces
in a suitable ordering.
Denote by $\chi_1,\dots,\chi_6:\fz_2^m\to\cz$ their characteristic
functions. There are the quadratic relations
$$\eqalign{
(\chi_1-\chi_2)(\chi_1-\chi_4)&=(\chi_3-\chi_6)(\chi_5-\chi_6)\cr
(\chi_1-\chi_2)(\chi_3-\chi_2)&=(\chi_5-\chi_4)(\chi_5-\chi_6)\cr
(\chi_1-\chi_2)(\chi_1-\chi_6)&=(\chi_3-\chi_4)(\chi_5-\chi_4)\cr
(\chi_1-\chi_2)(\chi_5-\chi_2)&=(\chi_3-\chi_4)(\chi_3-\chi_6)\cr
(\chi_3-\chi_4)(\chi_1-\chi_4)&=(\chi_5-\chi_2)(\chi_5-\chi_6)\cr
(\chi_3-\chi_4)(\chi_3-\chi_2)&=(\chi_1-\chi_6)(\chi_5-\chi_6)\cr
(\chi_1-\chi_4)(\chi_1-\chi_6)&=(\chi_3-\chi_2)(\chi_5-\chi_2)\cr
(\chi_1-\chi_4)(\chi_5-\chi_4)&=(\chi_3-\chi_2)(\chi_3-\chi_6)\cr
(\chi_1-\chi_6)(\chi_3-\chi_6)&=(\chi_5-\chi_2)(\chi_5-\chi_4)\cr
}$$
}
Qrd%
\finishproclaim
The proofs of \SsR\ and \Qrd\ are trivial.\qed
There are many other relations between characteristic functions.
But the mentioned here will improve to be important, because in some
cases they lead to relations between certain modular forms.
\smallskip
We remark that
the number of $(m-2)$-dimensional totally isotropic spaces in $\fz_2^{2m}$
is
$${5\over{2^m+1}} \cdot \prod_{i=1}^{m-2} {{(2^{2m+1-i}-1)}\over{(2^{m-1-i}-1)}}$$
Fortunately the 9 quadratic relations collaps to one if one takes
the linear relation \LRl\ into account. This gives us:
\proclaim
{Theorem}
{Let $A$ be a $(m-2)$-dimensional totally isotropic space.
Let $P$ be the five dimensional space which is generated by the
characteristic functions of maximal totally isotropic spaces containing
$A$ (see \LRl). There is a distinguished element $\pm R_A\in \Sym^2(P)$
which  for the standard $A$ can be read off from \Qrd\
and is defined for general $A$ in an obvious way
(using decompositions of a union of two stars as union of
two other stars).
This element is non-zero and contained in the kernel of
$$\Sym^2(P)\lo\cz[\fz_2^{2m}].$$}
TbK%
\finishproclaim
\zwischen{A homogenous ideal}%
For each maximal totally isotropic subspace $V\subset\fz_2^{2m}$ we consider a variable
$X_A$ and we consider the polynomial ring $\cz[\dots X_A\dots]$ in all
these variables.
We are interested in the subring
$\calR_m:=\cz[\dots X_A-X_B\dots]$.
which is generated by differences
$X_A-X_B$.
It is sufficient to take a fixed $A$ and arbitrary $B$. Hence
$\calR_m$ is a polynomial ring in one variable less.
For the generation of $\calR_m$ it is also sufficient
to take those with $ {\rm dim}(A\cap B)=m-1$.
Replacing the characteristic function by the corresponding
variable, we obtain from \LRl\ and \Qrd\
a finite set of linear and quadratic relations
of the form
$X_1+X_3+X_5-X_2-X_4-X_6$
(s.~\LRl) and from elements of the form
$(X_1-X_2)(X_1-X_4)-(X_3-X_6)(X_5-X_6)$ as explained in \Qrd.
They are contained in $\calR_m$ and we can consider the ideal
$$\calI_m\subset\calR_m$$
generated by them.
We believe that this ideal
is of interest and related to modular varieties. This paper contains
a contribution in the case $m=6$. There are related investigations in the
cases $m=3$ and $m=5$, cf. [Ko1] and [Ko2]. One knows
$$\dim(R_m/I_m)=\cases{1& if $m=1$\cr
3& if $m=2$\cr
6& if $m=3$\cr}$$
Unfortunately not much more is known about this ideal. A result which has
an application to modular varieties is that in the ideal certain quartic
relations can be found:
To explain this we consider an totally isotropic subspace $A\subset\fz_2^{2m}$ of
dimension $m-3$, where we assume $m\ge3$.
Let $B$ be   a coset which consists of anistropic elements
only. There are precisely 15 stars containing $B$. These 15 stars belong to
30 maximal totally isotropic spaces.
We describe them concretely. For this we can assume $\fz^{2m}=\fz^{2(m-3)}\oplus
\fz_2^6$ , $A$ is the standard totally isotropic space in $\fz^{2(m-3)}$
defined by $x_2=x_4=\cdots=x_{2m-6}=0$ and  $B$ is the coset  A+(0,\dots,0,1,1,0,0,0,0). Then the thirty spaces are of the form
$A\oplus X$, where $X$ is a maximal totally isotropic subspace in $\fz_2^6$. We list
the 30 spaces giving a basis in each case:
\vskip5mm\noindent
{\ninepoint\hbox to\hsize{\vbox{
\halign{\hfil#$\;$&$#\hfil$&$#\hfil$&$#\hfil$\cr
1&(0,0,0,0,0,1)& (0,0,0,1,0,0)& (0,1,0,0,0,0)\cr
2&(0,0,0,0,0,1)& (0,0,0,1,0,0)& (1,0,0,0,0,0)\cr
3&(0,0,0,0,0,1)& (0,0,1,0,0,0)& (0,1,0,0,0,0)\cr
4&(0,0,0,0,0,1)& (0,0,1,0,0,0)& (1,0,0,0,0,0)\cr
5&(0,0,0,0,0,1)& (0,1,0,1,0,0)& (1,0,1,0,0,0)\cr
6&(0,0,0,0,0,1)& (0,1,1,0,0,0)& (1,0,0,1,0,0)\cr
7&(0,0,0,0,1,0)& (0,0,0,1,0,0)& (0,1,0,0,0,0)\cr
8&(0,0,0,0,1,0)& (0,0,0,1,0,0)& (1,0,0,0,0,0)\cr
9&(0,0,0,0,1,0)& (0,0,1,0,0,0)& (0,1,0,0,0,0)\cr
10&(0,0,0,0,1,0)& (0,0,1,0,0,0)& (1,0,0,0,0,0)\cr
11&(0,0,0,0,1,0)& (0,1,0,1,0,0)& (1,0,1,0,0,0)\cr
12&(0,0,0,0,1,0)& (0,1,1,0,0,0)& (1,0,0,1,0,0)\cr
13&(0,0,0,1,0,0)& (0,1,0,0,0,1)& (1,0,0,0,1,0)\cr
14&(0,0,0,1,0,0)& (0,1,0,0,1,0)& (1,0,0,0,0,1)\cr
15&(0,0,0,1,0,1)& (0,0,1,0,1,0)& (0,1,0,0,0,0)\cr}}
\hfill\vbox{
\halign{#$\;$\hfil&#\hfil&#\hfil&#\hfil\cr
16&(0,0,0,1,0,1)& (0,0,1,0,1,0)& (1,0,0,0,0,0)\cr
17&(0,0,0,1,0,1)& (0,1,0,0,0,1)& (1,0,1,0,1,0)\cr
18&(0,0,0,1,0,1)& (0,1,1,0,1,0)& (1,0,0,0,0,1)\cr
19&(0,0,0,1,1,0)& (0,0,1,0,0,1)& (0,1,0,0,0,0)\cr
20&(0,0,0,1,1,0)& (0,0,1,0,0,1)& (1,0,0,0,0,0)\cr
21&(0,0,0,1,1,0)& (0,1,0,0,1,0)& (1,0,1,0,0,1)\cr
22&(0,0,0,1,1,0)& (0,1,1,0,0,1)& (1,0,0,0,1,0)\cr
23&(0,0,1,0,0,0)& (0,1,0,0,0,1)& (1,0,0,0,1,0)\cr
24&(0,0,1,0,0,0)& (0,1,0,0,1,0)& (1,0,0,0,0,1)\cr
25&(0,0,1,0,0,1)& (0,1,0,0,0,1)& (1,0,0,1,1,0)\cr
26&(0,0,1,0,0,1)& (0,1,0,1,1,0)& (1,0,0,0,0,1)\cr
27&(0,0,1,0,1,0)& (0,1,0,0,1,0)& (1,0,0,1,0,1)\cr
28&(0,0,1,0,1,0)& (0,1,0,1,0,1)& (1,0,0,0,1,0)\cr
29&(0,0,1,1,1,1)& (0,1,0,1,0,1)& (1,0,0,1,1,0)\cr
30&(0,0,1,1,1,1)& (0,1,0,1,1,0)& (1,0,0,1,0,1)\cr}}}
}

\proclaim
{Proposition}
{Let $X_i$ be the variables corresponding to the thirty maximal totally isotropic spaces
in $\fz_2^{2m}$ in the table above. The quartic polynomial
$$\eqalign{&
(X_{29}-X_{30})(X_{17}-X_{18})(X_{27}-X_{28})(X_{7}-X_{8})+\cr&
(X_{11}-X_{12})(X_{13}-X_{14})(X_{15}-X_{16})(X_{21}-X_{22})\cr}$$
is contained in the ideal $\calI_m$.}
QcQ%
\finishproclaim
We proceed as in \RbC\ , We decompose $\fz_2^{2m}=\fz_2^{2m-6}\oplus\fz_2^6$.
Let $C\in\fz_2^{2m-6}$
be a maximal totally isotropic subspace. Then $C$ can be considered as a
$ m-3$ dimensional totally isotropic space in $\fz_2^{2m}$. The maximal
totally isotropic subspaces which contain  it are of the form $C\oplus B$ where
$B\subset\fz_2^6$ is a three dimensional totally isotropic subspace of $\fz_2^6$.
So we reduced ourselves to the study of
totally isotropic subspaces in $\fz_2^6$ .
In this case we obtain $420$ quartic relations, cf [Ko2] ,
$36$ quadratic relations and $36$ linear relations
between the $105$ stars. Using the computer
algebra system SINGULAR,  we got that there are 14 linear independent
quadratic relations and the  generated ideal contains the quartics.
A similar calculation seems to be performed
in [Koi].
\neupara{Orthogonal modular forms}%
In the rest of this text, $L$ denotes an even lattice of signature $(2,n)$.
We consider the subset of $L\otimes_\gz\cz$ defined by
$$(z,z)=0,\quad (z,\bar z)>0.$$
It is an analytic manifold, which consists of two
connected components, which can be exchanged by the
map
$z\mapsto \bar z$. We denote by $\tilde\calH_n$
one of the components and by $\calH_n$ its image in the complex
projective space
$P(L\otimes_\gz\cz)$.
There is a subgroup of index two
$\O^+(V)$ of the orthogonal group $\O(V)$  ($V=L\otimes_\gz\rz$)
which preserves the two connected components.
This group acts on $\tilde H_n$ holomorphically.
For any $z=x+\imag y\subset\tilde\calH_n$ we consider the vector space
$$W=\rz x+\rz\imag y\in V:=L\otimes_\gz\rz.$$
This is  a positive definite subspace of $V$.
Moreover the space $W$ depends only on the image of $z$ in $\calH_n$.
This defines a bijection between $\calH_n$ and the set of two dimensional
positive definite subspaces of $V$.
\smallskip
The full group $\O(V)$ acts on the set of two-dimensional positive
definite subspaces. Using the bijection above, on
obtains an action of the whole $\O(V)$ (not only $\O^+(V)$)
on $\calH_n$.
But this action is non-holomorphic for
$g\in\O(V)-\O^+(V)$. Tt is given by
the real analytic map $z\mapsto \overline{g(z)}$.
\smallskip
We recall the notion of an
(holomorphic or meromorphic)
orthogonal modular form of weight $k\in\gz$ with
respect to a subgroup $\Gamma\subset\O^+(L)$ of finite index
and with respect to a character $v:\Gamma\to\cz^\bullet$.
It is a (holomorphic or meromorphic)
function $f:\tilde\calH_n\to \cz$ with
the properties
$$\eqalign{f(\gamma z)&=v(\gamma)f(z),\cr
  f(tz)&=t^{-k}f(z).\cr}$$
A  condition at the cusps has to be added. This condition is
in most cases automatically satisfied, for example for $n\ge 3$.
The vector space of holomorphic forms, which are
regular at the cusps, is denoted by $[\Gamma,k,v]$
and by $[\Gamma,k]$ when $v$ is trivial.
These are finite dimensional spaces
and moreover, the algebra
$$A(\Gamma)=\bigoplus_{k\in\gz}[\Gamma,k]$$
is finitely generated.
\smallskip
Let $\Mp(2,\gz)$  be the metaplectic cover of $\SL(2,\gz)$.
The elements of $\Mp(2,\gz)$ are pairs
$(M,\sqrt{c\tau+d}),$
where $M={ab\choose cd}\in\SL(2,\gz)$, and $\sqrt{c\tau+d})$
denotes a holomorphic root on  $\bfH$ of $c\tau+d$.
It is well known that $\Mp(2,\gz)$ is generated by
$$\eqalign{
T&= \left( {1\;1\choose 0\;1}, 1\right),\cr
S&= \left( {0-1\choose 1\phantom{-}0}, \sqrt{\tau}\right),\quad
\Re\sqrt\tau>0.}$$
One has the relations $S^2=(ST)^3=Z$, where
$Z=\bigl( {-1\phantom{-}0\choose \phantom{-}0-1}, \imag\bigr)$
is the standard generator of the center of $\Mp(2,\gz)$.
\smallskip
Recall that there is a unitary representation $\rho_L$
of $\Mp(2,\gz)$ on the group algebra $\cz[L'/L]$.
$$\eqalign{\varrho_L(T)&= \bigl(e^{2\pii q_L(\alpha)}
\bigr)_{\alpha\in L'/L}\quad
 \hbox{(diagonal matrix)}\cr
\varrho_L(S) &= { \sqrt{\imag}^{n -2}\over\sqrt{|L'/L|}}
\bigl(e^{-2\pii (\alpha,\beta)}\bigr)_{\alpha,\beta\in L'/L}.}$$
This representation is the Weil representation
attached to the finite quadratic form $(L'/L,\bar q_L)$.
\smallskip
We recall the notion of an elliptic modular form with respect to a
finite dimensional representation $\varrho:\Mp(2,\gz)\to\GL(W)$.
Let $k\in (1/2)\gz$ and $f:\bfH\to W$ be a
holomorphic function. Then $f$
is called {\it modular form\/} of weight $k$ with respect
to $\rho$ if
$$f(M\tau)=\sqrt{c\tau+d}^{2k} \varrho(M,\sqrt{c\tau+d})f(\tau)$$
for all $(M,\sqrt{c\tau+d})\in\Mp(2,\gz)$ and if
$f$ is holomorphic at $\imag\infty$.
We denote the space of all these modular forms by
$[\SL(2,\gz),k,\varrho]$.
A form is called a cusp form if it
vanishes at $\infty$. The space $[\SL(2,\gz),k,\varrho]$
decomposes into the direct sum of the subspace of cusp forms
and into the space of Eisenstein series which can be defined as
orthogonal complement with respect to the Petersson scalar product.
An Eisenstein series is determined by its constant Fourier coefficient.
\zwischen{The additive lift}%
In [Bo1],
Borcherds defined for integral $k+n/2$ a linear map
$$[\SL(2,\gz),k,\varrho_L]\lo[\Gamma_L,k+n/2-1],$$
which generalizes constructions of Saito-Kurokawa, Shimura,
Maa"s, Gritsenko, Oda et.al.
We are interested in this construction in the case $k=0$. Modular forms
of weight zero are constant, hence we get for even $n$ a map
$$\cz[L'/L]^{\SL(2,\gz)}\lo [\Gamma_L,n/2-1].$$
These orthogonal modular forms are the simplest examples of modular
forms of several variables. The weight $n/2-1$ is the so-called singular
weight. Every modular form of weight $0<k<n/2-1$ vanishes. Modular forms of singular
weight which are cusp-forms vanish identically%
\footnote{*)}{\ninepoint
The theory of singular modular forms is well established in the
Siegel case. For the orthogonal case there seems to be no good reference at
the moment. The proofs for the basic facts are the same.}.
\smallskip
There is a remarkable special case. Assume that $L$ is even and selfdual.
Then the Weil representation is the trivial one-dimensional representation
and we get a linear map $\cz\to [\O^+(L),n/2-1]$.
The image of 1 is a very distinguished modular form. Borcherd's formula
for the Fourier coefficients show that it is different from 0 and that
it agrees with a modular form which has already considered earlier by
Gritsenko and in a special case also by Krieg.
We call it here {\it Gritsenko's singular modular form.}
\smallskip
We recall  from [Fr2] some facts about the Baily Borel compactification
of $\calH_n/\Gamma$ for a subgroup of finite index $\Gamma\subset
\O^+(L)$. There are one- and zero-dimensional boundary components
which correspond to the $\Gamma$-orbits of two and one-dimensional
totally isotropic subspaces of $L\otimes_\gz\qz$.
We recall that the value of an orthogonal modular form
$F\in[\Gamma,k]$ at a non-zero isotropic element $\alpha\in
L\otimes_\gz\qz$ can be defined. One chooses $\beta\in L\otimes_\gz\qz$
such that $(\alpha,\beta)=1$
and defines
$$F(\alpha):\lim_{\Im\tau\to\infty}F(\tau\alpha+2\imag(\alpha,\beta)\beta).$$
It is easy to check that this limit exists (using the Fourier expansion
of $F$ as described in [Bo1]) and that this limit is independent of
the choice of $\beta$
and that it depends only on the $\Gamma$-orbit of
$\alpha$. One has
$F(t\alpha)=t^{-k}F(\alpha)$.
\smallskip
Representing the zero-dimensional boundary components by primitive
elements  $\alpha\in L'$ we obtain:
\smallskip\noindent
{\it When the weight $k$ is even, $F(\alpha)$ (for
$\alpha\in L'$ primitive isotropic) can be considered
as function on the set of zero-dimensional boundary points of
$\calH_n/\Gamma$.}
\smallskip\noindent
Additive lifts of constants have a basic property:
\proclaim
{Lemma}
{Let $F$ be a singular modular form, which is the additive lift
of a constant elliptic modular form. When all values
$F(\alpha)$, $\alpha\in L'$, primitive and isotropic, are zero, then
$F$ is identically zero.}
Lsz%
\finishproclaim
This follows from Borcherds' description of the Fourier
expansion of $F$ [Bo1], theorem 14.3. The formulae there
show that the "`restriction"'
of $F$ to the one-dimensional boundary components are
elliptic Eisenstein series, which hence are
orthogonal to cusp-forms. Eisenstein series vanish identically
if they vanish at all cusps. We don't want to give more details.
We only mention that Borcherds constructed examples of
non-trivial
singular modular forms which vanish at all zero dimensional cusps.
So these cannot be additive lifts. A systematic treatment of them is
a major unsolved problem.
\smallskip
Because the value at a cusp appears as constant Fourier coefficient,
we can use Borcherds formula from from [Bo1], theorem 14.3 to compute the values
of additive lifts at the cusps.
(In the first line of the final formula in 5. one has to replace
$c_{\delta z}(0)$ by $c_{\delta z/N}(0)$.)
We only reproduce the formula in a very special case.
Let $B(n)$ be the $n$-th Bernoulli  number.
\proclaim
{Proposition}
{Assume that $n\equiv 2\;\mod 4$ and that $2L'\subset L$. Let $F$
be the additive lift of $C\in\cz[L'/L]^{\SL(2,\gz)}$.
Then $F(\alpha)$ at a primitive
isotropic element $\alpha\in L'$ is
$$-{B(n/2-1)\over n-2}C(0)\qquad\hbox{if $\alpha\in L$}$$
 and otherwise
$$- {B(n/2-1)\over n-2}\bigl(C(0)+
C(\alpha)(1-2^{n/2-1})\bigr)
$$
}
KF%
\finishproclaim
We immediately obtain:
\proclaim
{Proposition}
{Assume that $n\equiv 2\;\mod 4$, $n>2$ and that $2L'\subset L$.
Assume furthermore that a linear subspace
$H\subset \cz[L'/L]^{\SL(2,\gz)}$ with the following property is given:
All elements of $H$ vanish at the zero element of $L'/L$.
Then there exists a non-zero constant $\gamma=\gamma(n)$ such that for
any $C\in H$ with corresponding additive
lift $F=F_C$ and every primitive isotropic $\alpha\in L'$
which is not contained in $L$ the formula
$$F(\alpha)=\gamma\cdot C(\alpha)$$
holds.
\smallskip\noindent
{\bf Corollary.} If in addition every non zero element
$\bar\alpha\in L'/L$ with $\bar q_L(\bar\alpha)=0$ can be represented
by an isotropic element of $L'$, then the additive lift
$$H\lo[\Gamma,n/2-1]$$
is injective.
}
Pfc%
\finishproclaim
Let $\alpha\in L'$ be an element of the dual lattice and $m<0$ a negative
integer. The Heegner divisor $H(\alpha,m)\subset\calH_n$ is the
union of all
$$v^\perp\cap\calH_n\qquad(v^\perp\hbox{
\rm orthogonal complement of $v$ in }P(V\otimes_\rz\cz)),$$
where $v$ runs through all elements in $L'$ with
$$v\equiv\alpha\,\mod\, L\quad\hbox{\rm and}\quad q_L(v)=m.$$
We consider $H(\alpha,m)$ as a divisor on $\calH_n$ by attaching
multiplicity $1$ to all components.
We have $H(\alpha,m)=H(-\alpha,m)$, more precisely, this divisor depends
only on the image of $\alpha$ in $(L'/L)/\pm1$.
\smallskip
The existence of an orthogonal modular forms whose zero divisor
is a given Heegner divisor is regulated by Borcherds'
{\it space of obstructions}
$$[\SL(2,\gz),n/2+1,\bar\varrho_L].$$
\proclaim
{Theorem}
{Assume $n>2$. Then the space of obstructions contains an
Eisenstein  series $\sum b_\alpha(m)e^{2\pii m\tau}$ with
constant Fourier coefficient
$$b_\alpha(0)=\cases{-1/2&if $\alpha=0$,\cr0&else.}$$
A finite linear combination
$$\sum_{\alpha\in (L'/L)/\pm1,\;m<0}
  C(\alpha,m)H(\alpha,m)\qquad(C(\alpha,m)\in\gz)$$
is the divisor of a meromorphic modular form $F$ of weight
$$k=  \sum_{m\in\gz,\,\alpha\in L'/L}b_\alpha(m)C(\alpha,-m)$$
if this number is integral and
if for every {\emph cusp form} $f$ in the the space of obstructions,
$$f_\alpha(\tau)=\sum_{m\in\qz} a_\alpha(m)\exp(2\pii m\tau),$$
the relation
$$\sum_{m<0,\,\alpha\in L'/L}a_\alpha(-m)C(\alpha,m)=0$$
holds.}
konteis%
\finishproclaim
(The assumption that $k$ is integral can be omitted if
one introduces modular forms of non-necessarily integral weight).
\neupara{  Level-two-cases}%
We want to consider level-two situations which means that we assume
$2L'\subset L$. Then
$L'/L$ is a vector space
over the field $\fz_2$ of two elements.
We use
the embedding
$$\fz_2\lo \qz/\gz,\qquad0\loma 0+\gz,\quad 1\loma 1/2+\gz.$$
The quadratic form $\bar q_L:L'/L\to\qz/\gz$ then can be considered
(as usual in the theory of quadratic forms over finite fileds)
as function
$$\bar q_L:L'/L\lo\fz_2.$$
For the rest of this paper we make the
\proclaim
{Assumption}
{The lattice $L$ is  even and has signature $(2,n)$. It has
the further property
$2L'\subset L$. The $\fz_2$-vector space $L'/L$ admits a basis such that
$$\bar q_L(x):=x_1x_2+\cdots+x_{2m-1}x_{2m}\qquad (2m=\dim(L'/L)).$$
(This means that the finite quadratic form is of ``even type''.)
}
Ann%
\finishproclaim
It is well-known that the signature $n-2$ mod 8 is determined by the
associated finite quadratic form. In our case it follows
$$n\equiv 2\ \mod 8.$$
The Weil representation on $\cz[\fz_2^m]$ factors through
$\SL(2,\gz/2\gz)$. It is a real representation with the  action of $\varrho_L (S)$
and
$$ \varrho_L (T) = \bigl((-1)^{ q(\alpha)}
\bigr)_{\alpha\in \fz_2^m}\quad
 \hbox{(diagonal matrix)}$$

We compute the character of this representation. Recall that
$\SL(2,\gz/2\gz)$ is isomorphic to the symmetric group $S_3$.
Hence it has three  conjugacy classes, which are characterized
by their orders 1,2,3.
Representatives are the matrices $E,S,ST$
(where $E$ denotes the unit matrix). The traces are given by
$$\eqalign{\trace(\varrho(E))&=2^{2m}\cr
 \trace(\varrho(T))&=2^{m-1}(2^m+1)-2^{m-1}(2^m-1)=2^m\cr
 \trace(\varrho(ST))&={1\over 2^m}(2^{m-1}(2^m+1)-2^{m-1}(2^m-1))=1\cr}$$
The character table is
$$\matrix{E&T&ST\cr\noalign{\vskip1mm}1&1&1\cr1&-1&1\cr2&0&-1\cr}$$
We obtain
\proclaim
{Lemma}
{We have
$$\dim\cz[\fz_2^{2m}]^{\SL(2,\gz)} = 2^{m-1}+{1\over 3}2^{2m-1}+{1\over 3}.$$
 In the case $m=6$ this is $715$.}
Lct%
\finishproclaim
On $\cz[\fz_2^{2m}]$ also  $\O(\fz_2^{2m})$
acts and this action commutes with the Weil representation.
Hence the group $\SL(2,\gz/2\gz)$ also acts
on the space of $\O(\fz_2^{2m})$-invariants
$$\cz[\fz_2^{2m}]^{\O(\fz_2^{2m})}.$$
This space is three dimensional and generated by three elements
$E_0,E_+,E_-$. Here $E_0$ means the characteristic function of the
zero element, $E_+$ is the characteristic function of the set of
non-zero isotropic elements and $E_-$ is the characteristic function
of the non isotropic elements.
Using this basis, we get a representation
$$\tilde\varrho:\SL(2,\gz/2\gz)\lo\GL(3,\cz).$$

From \Lai\ we obtain:
\item{1)} The matrix $\tilde\varrho(T)$ is the diagonal matrix with the
entries $1,1,-1$.
\item{2)} The diagonal entries of $\tilde\varrho(S)$ are
$${1\over 2^m},\quad {1\over 2}-{1\over 2^m},\quad {1\over 2}.$$
\vskip0pt\noindent
We obtain
$$\trace\tilde\varrho(E)=3,\quad \trace\tilde\varrho(T)=1,\quad
 \trace\tilde\varrho(ST)=0.$$
A glance at the character table shows:
\proclaim
{Lemma}
{The space
$$\cz[\fz_2^{2m}]^{\O(\fz_2^{2m})}$$
decomposes under the Weil representation into an irreducible
two dimensional and a trivial one-dimensional representation.
The space
$$\cz[\fz_2^{2m}]^{\O(\fz_2^{2m})\times\SL(2,\gz)}$$
is one dimensional. It is generated by the function which assigns
$(2^{m-1}+1)$ to the zero element, $1$ to all non-zero isotropic and $0$ to the
anisotropic elements of $\fz_2^{2m}$.}
Lgm%
\finishproclaim

\proclaim
{Proposition}
{The additive lift
$$\cz[\fz_2^{2m}]^{\SL(2,\gz)}\lo [\Gamma_{L},n/2-1]$$
is injective.}
Pii%
\finishproclaim
{\it Proof.\/} We decompose the space $\cz[\fz_2^{2m}]^{\SL(2,\gz)}$.
We consider the subspace $H$ of all functions $C$,
which satisfy the    linear equation
$$C(0)=0.$$
The full invariant element is not contained in $H$.
Hence $\cz[\fz_2^{2m}]^{\SL(2,\gz)}$ is the direct
sum of the one-dimensional invariant space and $H$ .
Since \Pfc\ applies, we obtain that the additive lift restricted to $H$
is injective.
The image contains no full invariant form. But the image of the invariant
element is a full invariant form, which is non-zero. Hence we have injectivitiy
on the whole.
\qed
\smallskipü
It is possible to construct explicit elements of $\cz[\fz_2^{12}]^{\SL(2,\gz)}$.
The following lemma is well-known, [Ko1], [Sch].
\proclaim
{Lemma}
{The  characteristic function of a maximal totally isotropic subspace of
$\fz_2^{2m}$
is $\SL(2,\gz)$-invariant.}
MIi%
\finishproclaim
The sum of all these characteristic functions is invariant under
the full modular group. Hence  we obtain:
\proclaim
{Lemma}
{
Assume that
$\cz[\fz_2^{2m}]^{\SL(2,\gz)}$ decomposes under $\O(\fz_2^m)$ into the sum
of a one-dimensional
and an irreducible subspace $H$.
Then
the space\/ $\cz[\fz_2^{2m}]^{\SL(2,\gz)}$ is spanned by the characteristic
functions of maximal totally isotropic subspaces.
\smallni
{\bf Remark.} The irreducibility of $H$ can be checked for $m\le 6$.
}
scS%
\finishproclaim

There is another important invariant element, which belongs to an
arbitrary non-zero isotropic element:
\proclaim
{Remark}
{Let $a$ be a non-zero isotropic element.
Let $\calA$ be the set of all maximal totally isotropic subspaces which contain
$a$ and $\calB$ the complementary set.
Denote by $\chi_A$ the characteristic
function of a maximal totally isotropic subspace. Consider the element
$$\sum_{A\in\calA}\chi_A-
 2^{m-2}\sum_{B\in\calB}\chi_B.$$
It is a multiple of the function that assumes the  value  $2^{m-2}$  on the
null  vector, $-2^{m-2}$  on the vector $a$, $1$ on the
isotropic vectors $b$ with $(a,b)=1$
and $0$ on the  remaining  vectors.}
Nvv%
\finishproclaim
We omit the easy proof. In the Special case $a=(1,0,\dots,0)$ we get
$$\cases{2^{m-2}& for $x=(0,0,\dots,0)$\cr
-2^{m-2}& for $x=(1,0,\dots,0)$\cr
1& for isotropic $x$ with $x_2=1$\cr
0&else\cr
}$$
\medskip\noindent
Now we consider the even and self dual lattice II$_{2,n}$, $n\equiv 2$ mod 8.
The discriminant group is zero.
To obtain a non-trivial discriminant group we rescale this lattice
and consider
$$L=\sqrt 2\;\hbox{\rm II}_{2,n},\qquad L'=L/2.$$
Then $L'/L$ is a vector space of dimension $2+n$ over the field
$\fz_2$ of two elements.
We can choose an isomorphism $\fz_2^{2+n}\to L'/L$  such
that   the quadratic form is
$$q(x):=x_1x_2+\cdots+x_{2m-1}x_{2m},\quad m=1+n/2.$$
Let us assume $n=10$, the group $\O(\fz^{12})$ contains a simple subgroup of index two.
The character table of this group is known.
For example the command\break
``CharTable(O12+(2).2)'' of the computer algebra system
GAP gives it.
The representations of
lowest dimensions have dimension $1$, $651$ and $714$ and they
are unique in these cases. Hence we obtain:
\proclaim
{Proposition}
{The additive lift of $\cz[\fz_2^{12}]^{\SL(2,\gz)}$ to $[\O^+(\hbox{\rm II}_{2,10}[2],4]$
is injective and
decomposes into the one-dimension trivial and the unique\/
$714$-dimensional irreducible  representation.}
Sev%
\finishproclaim
\neupara{Some Borcherds' products}%
We first investigate Borcherds' products for the full modular group
$\O^+(\hbox{\rm II}_{2,n})$, $n\equiv 2$ mod 8. The space of
obstructions is the space of usual elliptic modular forms
of weight $n/2+1$. The cases $n=10,18,26$  are of particular interest
because there are no elliptic cusp forms of weight 6, 10 and 14.
It is known that the group $\O^+(\hbox{\rm II}_{2,n})$
acts transitively on the set of all primitive vectors of a given norm.
We study the two lowest cases, vectors of norm $(x,x)=-2$ and
$(x,x)=-4$.
Before we continue, we modify slightly our notation for the
Heegner divisors. 
We denote by $H(m)$ the set of all
elements of $\calH_n$, which are orthogonal to some vector
$x\in \hbox{\rm II}_{2,n}$ with $(x,x)=m$.
For $\alpha\in\fz^{2+n}$ we denote by $H_\alpha(m)$ the set of
all points of $\calH_n$, which are orthogonal to some
$x\in \hbox{\rm II}_{2,n}$ with $(x,x)=m$ and $x\equiv\alpha$ mod 2.
Hence we have
$$H(m)=\bigcup_{\alpha\in\fz^{2+n}}H_\alpha(m).$$
With the notations of section 3 we have $H(m)=H(0,2m)$ where the
underlying lattice for the definition of $H(0,m)$ is
$\hbox{\rm II}_{2,n}$.
If we identify $\alpha$ with an element of $L'/L$ where
$L=\sqrt 2{\rm II}_{2,n}$ we get because of the rescaling factor
$$H_\alpha(m)=H(\alpha,2m).$$
The Heegner divisors
$H(-1)$ and $H(-2)$
are irreducible in $\calH_n/\O^+(\hbox{\rm II}_{2,n})$.
This follows from the fact that two primitive vectors of
${\rm II}_{2,n}$ of the same norm are equivalent
mod $\O^+(\hbox{\rm II}_{2,n})$ and that elements of norm $-2$
and $-4$ are automatically primitive.
We compute the weights of the modular forms with this divisor.
They are regulated by the Eisenstein series of weight
$6$, $10$ and $14$.
Their Fourier expansion, normalized such the constant
coefficient is $-1/2$, starts as follows:
$$\eqalign{\hbox{weight}\ 6:\hskip1.1mm\quad
&-{1\over 2}+252q+8316q^2+\cdots\cr
\hbox{weight}\ 10:\quad&-{1\over 2}+132q+67716q^2+\cdots\cr
\hbox{weight}\ 14:\quad&-{1\over 2}+12q+98316q^2+\cdots\cr}$$
\proclaim
{Lemma}
{In the cases $n=10,\;18,\;26$ there exists a
modular form on the full group $\O^+(\hbox{\rm II}_{2,n})$
whose zero divisor is the irreducible divisor $H(-1)$. The weights are
$252,\;132,\;12$. There also exists a  modular forms whose zero divisor
is the irreducible divisor $H(-2)$.
The weights are $8316,\;67716,\;98316$.}
Lmm%
\finishproclaim
We now investigate the level two case.
\proclaim
{Lemma}
{The divisor\/ $H_\alpha(-1)$ is not empty if and
only if $\alpha\in\fz_2^m$
is anisotropic.
In this case it is irreducible in
$\calH_n/\O^+(\hbox{\rm II}_{2,n})[2]$. Hence\/
$H(-1)$ considered in $\calH_n/\O^+(\hbox{\rm II}_{2,n})[2]$
has $2^{m-1}(2^m-1)$ irreducible components.
\smallskip\noindent
The divisor\/ $H_\alpha(-2)$ is not empty if
and only if $\alpha\in\fz_2^m$
is a non zero isotropic.
In this case it is irreducible in
$\calH_n/\O^+(\hbox{\rm II}_{2,n})[2]$. Hence\/
$H(-2)$ considered in $\calH_n/\O^+(\hbox{\rm II}_{2,n})[2]$
has $2^{m-1}(2^m+1)-1$ irreducible components.
}
Lhc%
\finishproclaim
We are interested in the case that the modular form, which belongs to
$H(-1)$ or $H(-2)$, splits into a
product of forms with irreducible
divisors $H_\alpha(-1)$ or $H_\alpha(-2)$.
But then the number of components
should divide the weight of this mo\-du\-lar form.
There is one distinguished case where this happens.
In the case II$_{2,10}$ and the divisor $H(-2)$ the weight has been
computed as 8316 (s.~\Lmm) and the number of components is 2079.
One has $8316/2079=4$. Thus we are lead to
\proclaim
{Proposition}
{For every non-zero isotropic $\alpha\in\fz_2^{12}$
(their number is $2079$)
there exists a modular form on the congruence group
of level two\/ $\O^+(\hbox{\rm II}_{2,10})[2]$
with divisor\/ $H_\alpha(-2)$.
The weight of this form is $4$.}
Pfs%
\finishproclaim
{\it Proof.\/}
We want to apply Borcherd's obstruction theory. For this purpose
we have to compute the space of elliptic cusp  forms of weight
6 with respect to the Weil representation on $\cz[\fz_2^{12}]$.
The Weil representation is trivial on the principal congruence
subgroup of level 2 of the elliptic modular group. The space
of cusp forms of weight 6 of this group has dimension one,
a generating element is $\eta^{12}$, where $\eta$ denotes the
Dedekind $\eta$-function. The expansion of $\eta^{12}$ is of the form
$$\eta(\tau)^{12}=\sum_{n>0} a_ne^{2\pii n\tau},\quad
a_n\ne 0\Longrightarrow n\in{1\over 2}+\gz.$$
Hence no $a_n\ne 0$ with integral $n$ occurs.
We see that actually all $H_\alpha(m)$ with {\it even\/} negative
$m$ (and non zero isotropic $\alpha$)
are divisors of modular forms.
\smallskip
We  shall consider particular additive lifts   of the elements from
$\cz[\fz_2^{12}]^{\SL(2,\gz)}$ and identify them with Borcherds products.
\proclaim
{Proposition}
{Let $M\subset\fz_2^{12}$ be a star (\Dsi).
The additive lift of the corresonding Weil invariant (\Kcs)  is a modular form
in $[\O^+(\hbox{II}_{2,10})[2], 4]$
whose zero divisor consists precisely of the $32$ Heegner divisors
$H_\alpha(-1)$ with $\alpha\in M$.}
Sta%
\finishproclaim
{\it Proof.\/}
Since the additive lift changes its sign under the  $2^{5}$
reflections, there must be zeros along the $2^{5}$ Heegner divisors.
We have to show that there are no other zeros. For this one takes the
product with respect to all $a$ and divides by a suitable power
with the form of weight 252 from \Lmm.
The result is a  form of weight 0 which has to be constant.\qed\smallskip
For similar statements cf [Ko1] and [Ko2].
\proclaim{Definition}
{Let $M$ be a star. The divisor consisting of the 32 Heegner divisors
$H_\alpha(-1)$ with $\alpha\in M$ is called the
associated {\emph star divisor}.}
DsD%
\finishproclaim
We consider now the additive lift of the Weil invariant, which has been
attached to an arbitrary non zero isotropic element in \Nvv. We will show that
this is also a Borcherds product.
\proclaim
{Theorem}
{Let $\alpha\in\fz_2^{12}$ be a non-zero isotropic vector.
The Borcherds product with divisor
$H_\alpha(-2)$ as described in \Pfs\ is contained in the
additive lift space. It is (up to a constant factor) the additive lift
of the Weil invariant defined in \Nvv.
}
Tab%
\finishproclaim
All what we have to show is that the additive lift vanishes along
$H_\alpha(-2)$. The proof of this is quite involved because
we have now reflection with fixes this divisor.
We postpone the proof to the next section.
\neupara{Projection of groups}%
Let $L$ be an even lattice of signature $(2,n)$ and $M\subset L$ a
sublattice of signature $(2,m)$. If the connected components
$\tilde\calH_n$, $\tilde\calH_m$ have been
properly chosen then there is a natural inclusion
$\tilde\calH_m\subset\tilde\calH_n$, which induces an inclusion
of the corresponding
half planes $\calH_m\subset\calH_n$.
Let $\Gamma\subset\O^+(L)$ a subgroup of finite index.
We consider the subgroup of all $g\in\Gamma$ with $g(M)=M$.
There is a natural homomorphism of this subgroup into $\O(M)$ and it is
easy to see that this subgroup ---let's denote it by $\Gamma'$---
is contained in $\O^+(M)$.
We call this group the {\it projected group.\/}
If $v,v'$ are compatible characters of $\Gamma,\Gamma'$
there is a natural restriction map
$$[\Gamma,k,v]\lo[\Gamma',k,v'].$$
\proclaim
{Theorem (Kneser)}
{Assume that $L$ is a lattice of index $(2,m)$, $m\ge 4$,
which has Witt-rank two
(which means that $L\times_\gz\qz$ contains a totally isotropic subspace of
dimension two). Assume furthermore that $L$ contains a vector of norm $-2$
and that the following two conditions hold:
\vskip2mm\item{\rm a)}
There exists a sublattice $L_1$ of rank $\ge 5$ whose discriminant
is not divisible by $3$.
\item{\rm b)} There exists a sublattice of rank $\ge 6$ with odd discriminant.
\smallni
Then the discriminant kernel $\Gamma_L$
is generated by reflections $x\mapsto x+(a,x)a$ along
norm $-2$ vectors $a$.
}
KR%
\finishproclaim
The lattice $D_m$ consists of all $x\in\qz^m$ such that
$x_1+\cdots+x_m$ is even. The quadratic form is
$$q(x)={1\over 2}(x_1^2+\cdots+x_m^2).$$
We apply Knesers result \KR\ to the lattice
$\gz^4\times(-D_m)$ where $m\ge 3$.
The lattices $A_1$ and $A_2$ can be embedded into $D(m)$, $m\ge 3$.
Hence we can take $L_1=\gz^4\times (-A_1)$
(discriminant 2) and $L_2=\gz^4\times (-A_2)$ (discriminant 3).
\proclaim
{Lemma}
{The discriminant kernels $\Gamma_L$ are generated by reflections along
vectors of norm $-2$ in the cases
$$L=L_m=\gz^4\times (-D_m)\quad\hbox{with}\quad m\ge 3.$$}
Ger%
\finishproclaim
We need another tool.\need4cm
\zwischen{A vanishing result}%
It may happen that a Heegner divisor
is the divisor of a Borcherds product.
In this case the weight of this product is given by a certain
Fourier coefficient of an elliptic Eisenstein series.
This Fourier coefficient is defined even if the Heegner divisor is not
a Borcherds product.
We denote this number the {\it virtual weight\/} of the Heegner divisor.
We refer to the paper [BK] of Bruinier and Kuss for the definition and
computation of these numbers.
We want to use the following non-trivial result of Bruinier:
\proclaim
{Proposition (Bruinier)}
{Let $L$ be an even lattice of signature $(2,m)$ with $m\ge 4$ and $f$
a non vanishing (entire) modular form of weight $k$
on the discrimant kernel $\Gamma_L$. Assume that $f$ vanishes along
a certain Heegner divsor $H$ (and maybe also on another other divisor).
Then $k$ is greater or equal than the
virtual weight of $H$.}
PB%
\finishproclaim

Now we have the tools for the
\smallni
{\it Proof of \Tab.}
We take the realization $\hbox{II}_{2,10}=\gz^4\times(-E_8)$
and we choose a concrete isomorphism
$\hbox{II}_{2,10}/2\hbox{II}_{2,10}\cong\fz_2^{12}$.
We consider a vector $a$ of norm -4. It can be taken inside $-E_8$.
It defines a non-zero isotropic vector $\alpha$ in $\fz_2^{12}$.
This vector defines a Weil invariant (s.~\Nvv). We denote by $f$ the additive lift
of it.
This is modular form with trivial character under the stabilizer
$\Gamma$ of $\alpha$ inside $\O^+(\hbox{II}_{2,10})$. From its concrete description
in connection with \Pfc\ follows that it vanishes at all cusps which are defined
by isotropic elements orthogonal to $a$.
\smallskip
All what we have to show is that $f$ vanishes along the orthogonal
complement of $a$. The orthogonal complement of a a norm 4 vector inside
$E_8$ is isomorphic to the root lattice $D_7$. Generally inside $D_{m+1}$
sits $D_m$ as orthogonal complement of a root.
Hence we can consider the chain
$$\gz^4\times(-D_2)\subset\cdots \gz^4\times(-D_7)\subset \gz^4\times(-E_8).$$
We denote the corresponding half planes by $\calH_4\subset\cdots\subset\calH_{10}$.
The claim is that $f$ vanishes on the nine dimensional $\calH_9$.
This will be done inductively starting form $\calH_4$.
\smallskip
The case $D_2$ is exceptional, one can not apply Kneser's theorem.
This lattice is generated by two orthogonal
vectors $(1,\pm1)$ of norm -2.
Hence it is isomorphic to $A_1\times A_1$.
From the explicit description
of the automorphism groups [CS] follows that
every automorphism of $D_2$
extends to an automorphism of $D_7$ and this extends
to an automorphism of $E_8$ which fixes $a$.
The orthogonal group $\O^+(\gz^4\times-( A_1\times A_1))$ is
isomorphic
to the extended Hermitian modular group of degree two
of the Gauss number field [FH].
This group is generated by Eichler transformations
and by the automorphisms of $A_1\times A_1$. Eichler transformations of course
extend
to $\O^+(\gz^4\times (-E_8)$. Hence we obtain that the projected group of
$\Gamma$ is the full $\O^+(\gz^4\times (-D_2))$ which can be identified
with the extended Hermitian modular group.
From the structure theorem [Fr3]
one can see every symmetric cusp form of weight 4 vanishes.
\smallskip
The next step is to prove that $f$ vanishes on $\calH_5$. It is a modular form
with respect to a certain projected group $\Gamma_5$. Now we can apply Kneser's result
to prove that this group contains the discrimant kernel of $\gz^4\times (-D_5)$.
The image of $\calH_4$ is contained in the Heegner divisor which is defined
by a root. But this Heegner divisor is irreducible. This follows from
the following well-known fact:
\proclaim
{Lemma}
{Let $L$ be an even lattice of which contains a copy of\/ \hbox{\rm II}$_{2,2}$.
Then any two primitive vectors $a,b\in L'$ with the same
image in $L'/L$ are equivalent mod $\Gamma_L$.}
TrE%
\finishproclaim
We have seen that the image of $\calH_4$ inside $\calH_5/\Gamma_5$ is a Heegner
divisor and that $f$ vanishes along it. Now we can apply Brunier's result \PB.
A concrete calculation with coefficients of Eisenstein series gives the virtual
weight $9$. Because this exceeds the actual weight 4, the form $f$ must vanish
on $\calH_5$. In the same way we prove the vanishing inductively on
$\calH_6,\dots,\calH_9$. The corresponding virtual weights are
$8$, $7$, $6$, $5$.
This completes the proof of \Tab.\qed
\neupara{Relations}%
In this  section we shall consider algebraic relations between  modular forms that are
additive liftings.
\smallskip
Since the additive lifting is a linear map, the linear relations
between the characteristic functions described
in \LRl\ of course imply the same relations between the corresponding
modular forms. But for the algebraic relations this is not clear.
We only can see that an algebraic relation between modular forms
implies the same relation between their values at the cusps. If we are inside
the space $H$ (defined by $C(0)=0$ then \Pfc\ shows that the input functions
$C\in (\fz_2^{2m})^{\SL(2,\gz)}$ (considered as functions on $\fz_2^m$) satisfy
the same relation. Hence there is a chance that the quadratic relations
in $(\fz_2^{2m})^{\SL(2,\gz)}$, which have been described in \Qrd\ also
hold for the corresponding modular forms. We don't know,
whether this is always true,
but it is at least true  when there
are modular forms whose zero set are star
divisors. This is the case when $n=10$ and $m=3,\ 5,\ 6$.
\bigskip
\proclaim
{Lemma}
{Let $A\subset\fz_2^{12}$ be a 4-dimensional totally isotropic subspace
above and $I_1,\dots,I_6$ the maximal totally isotropic subspaces
containing $A$ and in suitable ordering (as described in \Qrd.)
Denote by $f_1,\dots,f_6$ the additive lifts of their characteristic
functions. There are the quadratic relations
$$\eqalign{
(f_1-f_2)(f_1-f_4)&=(f_3-f_6)(f_5-f_6)\cr
(f_1-f_2)(f_3-f_2)&=(f_5-f_4)(f_5-f_6)\cr
(f_1-f_2)(f_1-f_6)&=(f_3-f_4)(f_5-f_4)\cr
(f_1-f_2)(f_5-f_2)&=(f_3-f_4)(f_3-f_6)\cr
(f_3-f_4)(f_1-f_4)&=(f_5-f_2)(f_5-f_6)\cr
(f_3-f_4)(f_3-f_2)&=(f_1-f_6)(f_5-f_6)\cr
(f_1-f_4)(f_1-f_6)&=(f_3-f_2)(f_5-f_2)\cr
(f_1-f_4)(f_5-f_4)&=(f_3-f_2)(f_3-f_6)\cr
(f_1-f_6)(f_3-f_6)&=(f_5-f_2)(f_5-f_4)\cr
}$$
These relations are induced by decompositions of a union of two
disjoint stars
as union of two other stars as for example $O_1\cup O_3=O_8\cup O_9$.
}
QrD%
\finishproclaim
{\it Proof.\/}
The point is that the divisors of both side agree by \Sta, because they are the
unions of the corresponding star divisors.
Hence both sides are equal up to a constant
factor. The normalizing factor follows from \Qrd.\qed
\smallskip
We recall that the 9 quadratic relations collaps to one if one takes
the linear relation \LRl\ into account.
\smallskip
In section 2 we introduced a certain ideal $\calI_{6}$
generated by linear and quadratic relations in the
polynomial ring
$\cz[\dots X_A-X_B\dots]$, where $X_A$ are formal
variables attached to maximal totally isotropic subspaces. We can summarize
\proclaim
{Proposition}
{There is a natural homomorphism
$$\calR_6/\calI_{6}\lo A(\O^+(\hbox{\rm II}_{2,10})[2]),$$
}
BaB%
\finishproclaim
\need 6cm
\neupara{Embedded lattices}%
Let us assume that the $L$ is  a sublattice of finite index
of another even lattice $M$.
Then we have
$$L\subset M\subset M'\subset L'.$$
The orthogonal groups $\O^+(M)$ and $\O^+(L)$ are usually not contained in each other
but for the discrimant groups one has
$$\Gamma_L\subset\Gamma_M.$$
To prove this one considers some element $g\in\Gamma_L$ and an element
$x\in M'$. Because $x$ is contained also in $L'$ we have $g(x)-x\in L\subset M$.
This shows $g(x)\in M'$. Hence $g$ defines an automorphism of $M'$ and then also
for $M=M''$. From $g(x)-x\in M$ we see $g\in \Gamma_M$.
\smallskip
Following Scheithauer we define an injective map
$$\Psi: \cz[M'/M]^{SL(2,\Bbb Z)} \to \cz[L'/L]^{SL(2,\Bbb Z)} .$$
as follows: We extend a function $M'/M\to\cz$ by zero to a function
$L'/M$ (i.e.~ the values are zero on the complement $L'/M-M'/M$)
and then compose the result with the natural projection
$L'/L\to L'/M$. Scheithauer [Sch] proved that Weil invariants go to Weil invariants.
\proclaim
{Lemma}
{The diagram
$$\matrix{\cz[M'/M]^{\SL(2,\gz)}&\lo &\cz[L'/L]^{\SL(2,\gz)}\cr
\downarrow&&\downarrow\cr
[\Gamma_M,n/2-1]&\subset&[\Gamma_L,n/2-1]\cr}$$
commutes.}
As%
\finishproclaim
{\it Proof.\/}
The proof follows from the definition of the theta lift [Bo1].
We explain it briefly for the lattice $M$. We consider a point $\tau$ form
the usual upper half plane and a $z\in\calH_n$. Let $A$ be the two dimensional positive
definit subspace of $V=M\otimes_\gz\rz$ and $B$ its orthogonal complement.
We decompose an arbitrary $v\in V$ as $v=v_A+v_B$ with $v_A\in A$ and $v_B\in B$.
The Siegel theta function corresponding to a coset $M+\alpha$ form $M'/M$ is
defined as
$$\theta_{M+\alpha}(\tau,z)=\sum_{\lambda\in M+\alpha}
e^{\displaystyle\pii\bigl(\tau(\lambda_A,\lambda_A)+\bar\tau(\lambda_B,\lambda_B)\bigr)}.$$
This can be considered as a function with values in $\cz[M'/M]$.
Using a standard pairing on $\cz[M'/M]$ we can pair the Siegel theta function
with an element $C\in\cz[M'/M]^{\SL(2,\gz)}$ and integrate then along
$\tau$ (with respect to an invariant volume element). The result essentially is
the additive lift of $C$. Now we compare the Siegel theta series with respect
to the two lattices $L\subset M$. For $\alpha\in M'$ one obviously has
$$\theta_{M+\alpha}=\sum_{\beta\in M/L}\theta_{L+\alpha+\beta}.$$
From this follows \As\ immediatly.
\qed
\proclaim
{Remark}
{Let $M$ be an even lattice and $L\subset M$ a sub-lattice of finite index.
The image of $\cz[M'/M]^{\SL(2,\gz)}$ under the map $\Psi$ consists of all
elements $C\in\cz[L'/L]^{\SL(2,\gz)}$, which are periodic under $M/L$.}
PIo%
\finishproclaim
{\it Proof.\/} Let $C\in\cz[L'/L]^{\SL(2,\gz)}$. Then
$$C(\alpha)={ \sqrt{\imag}^{n -2}\over\sqrt{|L'/L|}}\sum_{\beta\in L'/L}C(\beta)
e^{-2\pii (\alpha,\beta)}.$$
If $C$ is periodic we can write
$$C(\alpha)={ \sqrt{\imag}^{n -2}\over\sqrt{|L'/L|}}\sum_{\beta\in L'/M}C(\beta)
\sum_{\gamma\in M/L}e^{-2\pii (\alpha,\beta+\gamma)}.$$
The inner sum is zero if $\alpha$ is not contained in $M'$. Hence we can consider
$C$ as a function on $M'/M$ and it is clear that this function is invariant
under $\SL(2,\gz)$.\qed
\neupara{Distinguished points}%
The aim of this and next section is to prove the  generic
injectivity of the map from the modular variety
$\calH_{10}/\O^+(\hbox{\rm II}_{2,10})[2]$ to the projective
space of dimension $714$.
We recall some basic facts from [Ni] about primitive lattices.
\smallskip
Let $(L,(\cdot,\cdot))$ be an even lattice.
There is a natural isomorphism
$$L'\cong\Hom_\gz(L,\gz).$$
Let $M\subset L$
be a sublattice. Then there is a natural map $L'\to M'$.
The lattice $M$ is called {\it primitive\/} inside $L$ if $L/M$ is free.
Then the sequence
$$0\lo M\lo L\lo L/M\lo 0$$
splits and $L'\to M'$ is surjective. We mention that $M$ is primitive inside
$L$ if and only if
$$\qz M\cap L=M.$$
A lattice $M$ is called maximal if there is no even over-lattice
$M\subset\tilde M$, $M\ne \tilde M$, of the same rank.
We see that maximal lattices $M$ are primitive in every even over-lattice
$M\subset L$.
Lattices with square free determinant are maximal.
\smallskip
Recall that for a quadratic space of signature $(2,n)$, $n>0$, the group
$\O^+(V)$ is generated by reflections
along vectors of negative norm. Let $V=V_1\oplus V_2$ be an
orthogonal  decomposition into
a positive definite $V_1$ and a negative definite $V_2$. An automorphism
of $V$ which fixes this decomposition is contained in $\O^+(V)$
if and only of the determinant on $V_1$ is positive.\smallskip

Again we consider an even unimodular lattice $L$ and a primitive
sublattice $M\subset L$. We ask whether every primitive
sublattice $N\subset L$ which is isometric to $M$ can be obtained from
$M$ by applying an automorphism of $L$. For this we have to consider
the orthogonal complements $M^\perp$ and $N^\perp$ inside $L$. We know that there
discrimants agree. Hence they define the same genus. We now make the
rather strong assumption that $M^\perp$ is unique in its genus. Then
$M^\perp$ and $N^\perp$ are isometric lattices. We choose isometries
$$M\Isom N\quad\hbox{and}\quad M^\perp\Isom N^\perp.$$
They define an isometric isomorphism
$$M+ M^\perp\Isom N+N^\perp.$$
This isomorphism extends to $L$ if and only if the glue groups agree.
Recall that the glue group of $M+ M^\perp$ in $L$ is of the form
$$A=\set{(a,\sigma a),\quad a\in M'/M},$$
where $\sigma:M'/M\to (M')^\perp/M^\perp$
is an isomorphism (changing the sign of the quadratic form).
We assume that every isometric automorphism of $(M')^\perp/M^\perp$
extends to an element of $\O(M)$. Then we can obtain by a suitable
choice of the isomorphism $M^\perp\to N^\perp$ that the glue groups match.
We see:
\proclaim
{Remark (Nikulin)}
{Let $L$ be a unimodular lattice and $M\subset L$ a primitive sublattice
such that the lattice $M^\perp\subset L$ is unique in its genus.
Assume that every isometry of the discrimant group of $M^\perp$ extends
to an isometry of $M^\perp$.
Then the group
$\O(L)$ acts transitively on the set of isometric embeddings
$M\hookrightarrow L$.}
NIk%
\finishproclaim
This applies to the following situation. Let $L=\hbox{\rm II}_{2,10}$
and $M=A_2$. We can take the realization $L=\hbox{\rm II}_{2,2}\times (-E_8)$
and realize $A_2$ inside $\hbox{\rm II}_{2,2}$, for example as the lattice
generated by $(1,1,0,0;0,\dots,0)$ and $(1,0,1,1;0,\dots,0)$.
We see that a realization of $A_2$ in $L$ exists such that
the orthogonal complement is $-(A_2\times E_8)$. From Kneser's classification
of definite lattices of determinant $\le 3$ up to dimension 17 follows
that this lattice is unique in its genus.
The discriminant is isomorphic to $\gz/3\gz$. The non trivial automorphism
$x\mapsto -x$ extends.
Hence \NIk\ applies. We can even more replace $\O(L)$ by $\O^+(L)$ because
$A_2$ admits an automorphism of determinant -1.
\proclaim
{Lemma}
{There is precisely one $\O^+(L)$-%
orbit of sublattices $M\subset
L=\hbox{\rm II}_{2,10}$ isometric to $A_2$.
}
Aze%
\finishproclaim
We want to determine the subgroup of $\O^+(L)$ which fixes a given $A_2$.
We start with the full group $\O(L)$. This group acts on the set
of two dimensional positive definite subspaces in a natural way.
The automorphism group of $A_2$ has order 12.
By a result of Eichler one knows
$$\O(A_2\times E_8)\cong \O(A_2)\times \O(E_8).$$
The order of $\O(E_8)$ is $696\,729\,600$.
We get
$$\#\O(A_2\times(-(A_2\times E_8)))=12\cdot12\cdot696\,729\,600= 2^{18}\cdot3^7\cdot5^2\cdot7.$$
Only half of this group survives in $\O(L)$ because of the matching condition.
\proclaim
{Lemma}
{The subgroup of $\O(\hbox{\rm II}_{2,10})$, which stabilizes a given
$A_2$, has order $2^{17}\cdot3^6\cdot5^2\cdot7$.}
SAz%
\finishproclaim
Next we consider the congruence group
$$\O(\hbox{\rm II}_{2,10})[2]:=\kernel\bigl(\O(\hbox{\rm II}_{2,10})
\lo\O(\fz_2^{12})\bigr).$$
The number of all $\O(\hbox{\rm II}_{2,10})[2]$-orbits of $A_2$-sublattices
is the index
$[\O(\fz_2^{12}):H]$, where $H$ denotes the image of the stabilizer
described in \SAz. To get the order of $H$ we have to know the order
of the kernel of the homomorphism
$$\O(A_2\times(-(A_2\times E_8)))\to \O(\fz_2^{12}).$$
When $M$ is a primitive sublattice of an even lattice $N$, one has
$M\cap 2N=2M$. This shows that in the kernel there are only signs changes
of the occurring $A_2$, $-A_2$ and $-E_8$. Hence the kernel has
order 8. Because of the gluing condition only order 4 survives.
From \SAz\ we obtain
$\#H=2^{15}\cdot3^6\cdot 5^2\cdot 7$.
Hence the number of $\O(\hbox{\rm II}_{2,10})[2]$-orbits of $A_2$-sublattices
is
$2^{16}\cdot 3\cdot 7\cdot 11\cdot 17\cdot 31$.
Let $M$ be a lattice which is isomorphic to $A_2\times E_8$, then there is a
unique sublattice of $M$ which is ismorphic to $A_2$ and whose
orthogonal complement is isomorphic to $E_8$.
This follows from the mentioned
result of Eichler. Hence an $A_2$-lattice inside $\hbox{\rm II}_{2,10}$ defines
a unique sublattice isomorphic to $-A_2$. There images in $\fz_2^{12}$
define an ordered pair $(A,B)$ of orthogonal hyperbolic planes.
Every pair occurs because the orthogonal group acts transitively on them.
The number of these pairs is easy to compute, one obtains
$2^{16}\cdot 3\cdot 7\cdot 11\cdot 17\cdot 31$, the same number as above.
Hence we obtain:
\proclaim
{Lemma}
{The\/  $\O(\hbox{\rm II}_{2,10})[2]$-orbits of sublattices of\/
$\hbox{\rm II}_{2,10}$ isomorphic to $A_2$ are in one-to-one correspondence
to the ordered pairs $(A,B)$ of orthogonal hyperbolic planes
in $\fz_2^{12}$. Their number is $2^{16}\cdot 3\cdot 7\cdot 11\cdot 17\cdot 31$.}
Smz%
\finishproclaim
The situation for $\O^+(\hbox{\rm II}_{2,10})[2]$ is slightly different.
The reason is that an automorphism which stabilizes an embedded $A_2$
is in $\O^+$ if and only if it has positive determinant on $A_2$.
Hence the number of $\O^+(\hbox{\rm II}_{2,10})[2]$-orbits is the double
of the number given in \Smz. This can be understood as follows:
\proclaim
{Definition}
{Two $\O^+(\hbox{\rm II}_{2,10})[2]$-classes  of sublattices
of $\hbox{\rm II}_{2,10}$ of type $A_2$ are called {\emph companions},
if there exists an element $g\in \O(\hbox{\rm II}_{2,10})[2]$,
which is not contained
in $\O^+(\hbox{\rm II}_{2,10})[2]$ such that $g$ maps one class to the other.}
DCc%
\finishproclaim
As an example we realize
$\hbox{\rm II}_{2,10}$ as $\hbox{\rm II}_{2,2}\times(-E_8)$.
We consider the lattice generated by $(1,1,0,0;0,\dots,0)$ and
$(1,0,1,1;0,\dots,0)$. The automorphism which changes the signs of the first two
components is not in $\O^+$. Hence the companion class is represented
by the lattice with basis $(-1,-1,0,0;0,\dots,0)$ and
$(-1,0,1,1;0,\dots,0)$.
\proclaim
{Remark}
{Companions are equivalent mod $\O(\hbox{\rm II}_{2,10})[2]$
and also equivalent mod $\O^+(\hbox{\rm II}_{2,10})$ but not equivalent mod
$\O^+(\hbox{\rm II}_{2,10})[2]$.}
Ecc%
\finishproclaim
Our method also shows that the stabilizer of a lattice of type $A_2$
(embedded in $\hbox{\rm II}_{2,10}$) inside $\O^+(\hbox{\rm II}_{2,10})[2]$
only consists of the identity and its negative.
We also can restate our results in terms of modular varieties:
\proclaim
{Proposition}
{The set of points of type $A_2$ in
$$\calH_{10}/\Gamma[2]\qquad(\Gamma[2]=\O^+(\hbox{\rm II}_{2,10})[2])$$
consists of\/ $2^{16}\cdot 3\cdot 7\cdot 11\cdot 17\cdot 31$ pairs of
companions. The map $\calH_{10}\to \calH_{10}/\Gamma[2]$ is
locally biholomorphic at
the $A_2$-points. Especially they define smooth points
in $\calH_{10}/\Gamma[2]$.}
SoP%
\finishproclaim
\neupara{Mapping to the projective space}%
We have to consider the Baily-Borel compactification
$\overline{\calH_{10}/\Gamma[2]}$, which is an algebraic variety [BB].
We want to use a basis $F_0,\dots,F_{713}$ of the
714-dimensional irreducible part of the additive lift space to define
a regular map from it to $P^{713}(\cz)$. For this we need
\proclaim
{Lemma}
{The modular forms of the 714-dimensional space of additive lifts
don't have common zeros in $\overline{\calH_{10}/\Gamma[2]}$.}
Kgn%
\finishproclaim
{\it Proof.\/}
The more involved part is the interior. Here
the argument is the same as in Kondo's case and we omit it.
(Actually it can be obtained as a consequence of Kondo's case
using \KU.) For the boundary one can argue as follows. The irreducible
components of the boundary are modular curves corresponding
to $\SL(2,\gz)[2]$. The restriction of the 714-dimensional space
contains a form of weight 4, which is not invariant under the full modular group.
The space of all elliptic modular forms is three-dimensional and generated
by
$\vartheta_1^8$, $\vartheta_2^8$, $\vartheta_3^8$, where
$\vartheta_1,\vartheta_2,\vartheta_3$ denote the Jacobi theta constants
with respect to the characteristics $(0,0)$, $(1,0)$, $(0,1)$. Recall that
$\vartheta_1^4=\vartheta_2^4+\vartheta_3^4$.
The space generated by the the three eighth powers splits  under $\SL(2,\gz)$
into a one dimensional space generated by
$f=\vartheta_1^8+\vartheta_2^8+\vartheta_3^8$ and an irreducible two dimensional
space generated by $\vartheta_1^8-\vartheta_2^8$,
$\vartheta_1^8-\vartheta_3^8$ and $\vartheta_2^8-\vartheta_3^8$. This space must
be in the restriction of the 714-dimensional space. But this two dimensional space
has no joint zero. This follows from the fact that $f^2$ is in the second symmetric
power of this two dimensional space (consider the sum of the
$(\vartheta_i^8-\vartheta_j^8)^2$).
\qed\smallskip
Due to a result of Hilbert we obtain from \Kgn:
The map
$$\overline{\calH_{10}/\Gamma[2]}\lo P^{713}(\cz)$$
defined by a basis of the additive lift space is a finite regular map.
We denote by
$$A(\Gamma[2])=\sum_{r=0}^\infty[\Gamma[2],r]$$
the graded algebra of modular forms (similarly for other groups) and
by
$$B(\Gamma[2])=\cz[F_0,\dots, F_{714}]\subset A(\Gamma[2])$$
the subring which is generated by the 715-dimensional
additive lift space. Be aware that we include the invariant form and not only
the 714-dimensional part.
\proclaim
{Lemma}
{The finite map
$$\overline{\calH_{10}/\Gamma[2]}=\proj(A(\Gamma[2]))\lo
\proj(B(\Gamma[2]))\qquad(\hookrightarrow
P^{714}(\cz))$$
is locally biholomorphic at the $A_2$-points.}
nSa%
\finishproclaim
{\it Proof.\/}
Let $A_2\times(-A_2\times E_8)\hookrightarrow \hbox{\rm II}_{2,10}$
an embedding which describes an $A_2$-point. In $-(A_2\times E_8)$ exist
10 linearly independent vectors of norm -4. We know that
there exist ten modular forms in the additive lift space whose zero
divisor consists of the corresponding Heegner divisors. These define locally
at the point a coordinate frame.\qed

We want to proof that the map in \nSa\ is generically injective.
We follow closely Kondo ([Ko1], section 6).
First of all we consider for a two-dimensional positive definite subspace
$W\subset\hbox{\rm II}_{2,10}\otimes_\gz\rz$ the set of all norm -2 vectors
in $\hbox{\rm II}_{2,10}$, which are orthogonal to $W$. The image of this
set in $\fz_2^{12}$ is denoted by
$$\Delta_W:=\set{\alpha\in\fz_2^{12};\quad\alpha\
\hbox{image of a norm -2 vector}\ a\in\hbox{\rm II}_{2,10}\
\hbox{orthogonal to}\ W}.$$
If $W$ is an $A_2$-point (i.e.~the vector space generated by a $A_2$-sublattice
of $\hbox{\rm II}_{2,10}$), the set $\Delta_W$ contains 123 elements.
This can be seen if one takes the realization
$$\hbox{\rm II}_{2,10}=\hbox{\rm II}_{2,2}\times (-E_8)$$
and for $W$ the space generated by $(1,1,0,0;0,\dots,0)$ and
$(1,0,1,1;0,\dots,0)$. The orthogonal norm -2 vectors are
Vectors  of norm -2 orthogonal to them are
$\pm(0,0,-1,1;0)$,
$\pm(1,-1,-1,0;0)$,
$\pm(1,-1,0,-1;0)$ and $(0,0,0,0,x)$
with $x$ a root of $- E_8$. This description also shows:
\proclaim
{Remark}
{Let $\varrho$ be an $A_2$-point and $A,B\subset\fz_2^{12}$ the corresponding
pair (\Smz). The set $\Delta_\varrho$ consists of two types:%
\smallskip\noindent
{\emph First type:\/} Three points inside $A\oplus B$.\hfill\break
{\emph Second type:\/} 120 points orthogonal to $A\oplus B$.}
Fst%
\finishproclaim
The set $\Delta_W$ only depends on the $\Gamma[2]$-equivalence class of $W$.
Hence we can define the set $\Delta_\varrho$ also for points
$\varrho\in\calH_{10}/\Gamma[2]$:
\proclaim
{Lemma}
{Let $\varrho,\varrho'\in\calH_{10}/\Gamma[2]$ be two points. Assume that one
of them is an $A_2$-points and furthermore
that $\Delta_\varrho=\Delta_{\varrho'}$.
Then the two points are equal or companions.}
cDA%
\finishproclaim
The following two lemmas are similar to Kondo's 6.4 and 6.5 in [Ko1]:
\proclaim
{Lemma}
{Let $\Delta=\Delta_\varrho$ for a $A_2$-point $\varrho$. Assume that
$\alpha\in\fz_2^{12}$ is a non-isotropic element, which is not contained
in $\Delta$. Then there exists a star, which contains $\alpha$ and
has empty intersection with $\Delta$.}
LKa%
\finishproclaim
This is a finite problem which can be checked by calculation.
One can assume that $\varrho$ is our standard point. There is
a group isomorphic to $S_3\times\O(E_8/2E_8)$ which stabilizes
$\Delta$ and also its complement. The complement consists of 5
orbits, which can be represented by
\vskip2mm
\halign{\quad#\hfil\quad&#\hfil\cr
$(0,0,1,1;0,1,0,0,0,0,0,0)$ &$(1,1,0,0;0,0,0,0,0,0,0,0)$\cr
$(1,1,0,0;0,0,1,0,0,0,0,0)$& $(0,0,0,1;1,1,0,0,0,0,0,0)$\cr
$(0,1,0,0;1,1,0,0,0,0,0,0)$&\cr}
\smallni
In the first four cases we take the following star:
The underlying totally isotropic space is spanned by
\vskip2mm\halign{\quad#\hfil\quad&#\hfil\cr
$(1,1,1,1;0,1,0,0,0,0,0,0)$&$(1,1,0,1;1,1,0,0,0,0,0,0)$\cr
$(0,0,0,0,0,0,1,0,0,0,0,0)$&$(0,0,0,0,0,0,0,0,1,0,0,0)$\cr
$(0,0,0,0,0,0,0,0,0,0,1,0)$&\cr}
\smallni
The star is the coset which contains $(1,1,0,\dots,0)$.
In the last case, we consider   the even totally isotropic space
spanned by
\vskip2mm\halign{\quad#\hfil\quad&#\hfil\cr
$(1,1,1,1,0,1,0,0,0,0,0,0)$&$(0,1,1,1,1,1,0,0,0,0,0,0)$\cr
$(0,0,0,0,0,0,1,0,0,0,0,0)$&$(0,0,0,0,0,0,0,0,1,0,0,0)$\cr
$(0,0,0,0,0,0,0,0,0,0,1,0)$&\cr}
\smallni
The star is the coset which contains
$(0,0,1,1;0,0,\dots,0)$.\qed
\proclaim
{Lemma}
{Let $\Delta=\Delta_\varrho$ for a $A_2$-point $\varrho$. Assume that
$\alpha\in\Delta$ is a non-isotropic element.
\vskip1mm\item{\rm a)} Assume that $\alpha$ is of the first type (see \Fst).
There exists a star $S$ such that
$S\cap\Delta=\{\alpha\}$.
\item{\rm b)} Assume that $\alpha,\beta\in\Delta$ are two orthogonal
elements of the second type.
There
exists a star $S$ such that $S\cap\Delta=\{\alpha,\beta\}$.
\vskip0pt}
LKb%
\finishproclaim
{\it Proof.\/} Again we assume that $\varrho$ is the standard point.
Let $\alpha$ be of first type.
Up to the action of the group $S_3$ we can take
$\alpha=(0,0,1,1;0,\dots,0)$.
We can take the star
$$(*,0,1,1,*,0,*,0,*,0,*,0).$$
Let $\alpha,\beta\in\Delta$ be orthogonal elements of second type.
Using the action of
$\O(E_8/2E_8)$, we can assume
$$\alpha=(0,0,0,0;1,1,0,\dots,0),\quad \beta=(0,0,0,0;1,1,0,0,0,0,1,0).$$
In this case we can take the star
$\alpha+M$ where $M$ is the totally isotropic space spanned by
the vectors
\halign{\quad#\hfil&\quad#\hfil&\qquad\hfill#\cr
$(1,0,0,0;0,0,1,0,0,0,0,0),$& $(0,1,0,0;0,0,0,1,0,0,0,0),$\cr
$(0,0,1,0;0,0,0,0,1,0,0,0),$& $(0,0,0,1;0,0,0,0,0,1,0,0,),$\cr
$(0,0,0,0;0,0,0,0,0,0,0,1).$&&\hfill\square\cr}
\smallni
Before we continue we recall that a non-isotropic element
$\alpha\in\fz_2^{12}$ defines a certain
Heegner divisor $H_\alpha(-1)$ in $\calH_{10}/\Gamma[2]$.
It is the fixed point set of the reflection
corresponding to a norm -2 representative
$a\in\hbox{\rm II}_{2,10}$. Hence we call this divisor a {\it mirror}.\/
We also notice that a point $\varrho$ is contained
in this mirror if and only if $\alpha\in\Delta_\varrho$.
If $S\subset\fz_2^{12}$ is a star, we can consider
the  star divisor which consists of the corresponding 32 mirrors.

The star divisor of a star $S\subset\fz_2^{12}$
contains $\varrho$ if and only of $S\cap\Delta_\varrho$ is not empty.
\proclaim
{Lemma}
{Let $\varrho\in\calH_{10}/\Gamma[2]$ be an $A_2$-point and let
$\varrho'$ be an arbitrary point which is different from $\varrho$ and its
companion.
Then there exists a star which contains one of the two points but not
both.}
Sdc%
\finishproclaim
{\it Proof.\/} We assume that
the two points cannot be separated.
If there exists an
$\alpha\in \Delta_{\varrho'}-\Delta_\varrho$ we find by \LKa\ a star
which contains $\varrho'$ but not $\varrho$. This shows
$\Delta_{\varrho'}\subset \Delta_{\varrho}$.
Let be
$\calM=\Delta_\varrho-\Delta_{\varrho'}$.
Every $\alpha\in\calM$ must be of
second kind, because otherwise we could use \LKb, a) to separate the
two points.
In $\calM$ two elements cannot be orthogonal, because otherwise we
could separate the two points using \LKb, b).
The same argument shows that
every  element of second type
which is orthogonal to an element of $\calM$
is contained in $\Delta_{\varrho'}$.
Before we continue we fix a notation:
\smallskip
An $E_8$-root system in $E_8/2E_8$ is a system of
anisotropic vectors
$\alpha_1,\dots,\alpha_8$
with the property that
$$(\alpha_1,\alpha_2)=(\alpha_2,\alpha_3)=\cdots(\alpha_7,\alpha_8)=1
\quad\hbox{and}\quad (\alpha_5,\alpha_8)=1.$$
{\it Claim.\/} Let $\calM\subset E_8/2E_8$ be a subset containing
only anisotropic elements, such that
any two elements of $\calM$ are not orthogonal.
Then the complement of $\calM$ contains an  $E_8$-root system.
\smallni
{\it Proof of the claim.\/}
The cases where $\calM$ is empty or contains only one element is trivial.
Hence $\calM$ contains at least two elements $\alpha,\beta$.
We use the standard realization $E_8/2E_8=\fz_2^8$ with
$q(x)=x_1x_2+\cdots+x_7x_8$
and can then assume
$$\alpha=(1,1,0,0,0,0,0,0,),\quad\beta=(0,1,1,0,1,0,1,1).$$
The system\vskip1mm
\halign{\ $#$\hfil&\ $#$\hfil&\ $#$\hfil&\ $#$\hfil\cr
(1,1,1,0,0,0,0,0)&(0,0,1,1,1,0,0,0)&(0,0,0,0,1,1,1,0)&(0,0,0,0,0,0,1,1)\cr
(1,1,0,0,0,0,1,0)&(0,1,0,1,1,1,0,0)&(0,1,0,1,0,1,1,1&(1,0,1,1,0,0,0,0)\cr}%
\smallni
is an $E_8$-root system. The first five are orthogonal to $\alpha$ the rest
to $\beta$. Hence they are all in the complement.\qed
\smallni
{\it Proof of \Sdc\ continued.\/}
We want to apply this to the following situation. As above
a representative of $\varrho$ belongs to a decomposition
$\hbox{\rm II}_{2,10}=\hbox{\rm II}_{2,2}\times (-E_8)$
The set $\calM=\Delta_\varrho-\Delta_{\varrho'}$ belongs to
the part $E_8/2E_8$ in this decomposition (second type).
From the above follows that $\Delta_{\varrho'}$ contains an $E_8$-root
system. Together with the three elements of the first type we obtain:
The set $\Delta_{\varrho'}$ contains an $A_2\times E_8$-root system
in an obvious sense. Now we have to recall the definition of
$\Delta_{\varrho'}$. It is the image of {\it all\/}  elements of norm -2
contained in a ten-dimensional negative definite
sublattice of $\hbox{\rm II}_{2,10}$. From the classification of
root lattices by their Dynkin diagram we see that the lattice generated by
these roots is a copy of $-(A_2\times E_8)$. This means that $\varrho'$
is an $A_2$-point. Hence $\Delta_{\varrho'}$ contains 123 elements
and we get $\Delta_\varrho=\Delta_{\varrho'}$. Now \cDA\
completes the proof of \Sdc.\qed\smallskip
We also have to take boundary points into account. Recall that there
are
one-dimensional and zero dimensional boundary components. Each one-dimensional
boundary component corresponds to a two dimensional totally isotropic space $B$.
It belongs to the group $\SL(2,\gz)[2]$ and the three cusp classes of this
group correspond to the three non-zero isotropics of $B$. Let
$\alpha\in\fz_2^{12}$ be an anisotropic element, which is orthogonal to
$B$. Then the three zero-dimensional boundary points are contained in the closure
of the Heegner divisor $H_\alpha(-1)$.
The number of these $\alpha$ is 480. Let $\varrho$ now be any $A_2$-point. As we know
it is contained in 123 Heegner divisors $H_\beta(-1)$. Hence there exists an
$\alpha$ such that the closure of $H_\alpha(-1)$ contains the three cusps but not
$\varrho$. Now we use \LKa\ and obtain a star modular form $f_S$ with the
following property: It doesn't vanish at $\varrho$ but it vanishes in the three
cusps in consideration. Now we use that every cusp from
of weight 4 for the group  $\SL(2,\gz)[2]$ vanishes. Hence we obtain:
\proclaim
{Lemma}
{Let $\varrho$ be an $A_2$-point. For every one-dimensional boundary
component there exist a star modular form which vanishes along this
boundary component but does not vanish at $\varrho$.}
TBk%
\finishproclaim
From \Sdc, \TBk\ and \nSa\ we know that the map
$\overline{\calH_{10}/\Gamma[2]}\to\proj(B(\Gamma[2]))$ has covering
degree at most two. We want to show that it is one.
Because we cannot separate companions directly by our methods, we have to take a
small detour.
We take invariants under the group $\O(\fz_2^{12})$ to go down
to the full modular group $\Gamma=\O^+(\hbox{\rm II}_{2,10})$.
We claim the this finite group acts faithfully on $B(\Gamma[2])$. This follows from
the fact that $\O(\fz_2^{12})$ has a simple subgroup of index 2 (the kernel of the
so-called Dixon invariant) and that this subgroup cannot act trivially on
$B(\Gamma[2])$. Set
$$B(\Gamma):=B(\Gamma[2])^{\O(\fz_2^{12})}.$$
We consider the diagram
$$\matrix{
\overline{\calH_{10}/\Gamma[2]}&\lo&\proj(B(\Gamma[2])\cr
\downarrow&&\downarrow\cr
\overline{\calH_{10}/\Gamma}&\lo&\proj(B(\Gamma)
\cr}$$
Since both vertical arrows have the same covering degree, it is sufficient to show
that the second row has covering degree one. But this follows from
\proclaim
{Lemma}
{
The finite map
$$\overline{\calH_{10}/\Gamma}=\proj(A(\Gamma))\lo
\proj(B(\Gamma))$$
is locally biholomorphic at the $A_2$-point. The $A_2$-point is unique in its fibre.}
NSa%
\finishproclaim
This follows from our considerations about the level two case and from
the fact that companions have the same image in level 1.
\proclaim
{Theorem}
{The map
$$\overline{\calH_{10}/\Gamma[2]}\to\proj(B(\Gamma[2]))$$
is everywhere regular finite and birational. The ring of modular forms
of weight divisible by 4
$$A^{(4)}(\Gamma[2]):=\bigoplus_{r=0}^\infty[\Gamma[2],4r]$$
is the normalization of the ring
$B(\Gamma[2])$ (the ring generated by the 715 dimensional additive
lift space). Using \BaB\ we get
a finite map
$$\calH_{10}/\O^+(\hbox{\rm II}_{2,10})[2]\lo
\proj(\calR_6/\calI_{6}).$$
}
RFb%
\finishproclaim

\neupara{Enriques surfaces, Kondo's approach}%
Denote by $U$ the unimodular lattice $\gz\times\gz$ with quadratic form
$(x,x)=2x_1x_2$.
Kondo investigated in [Ko1] the case of the lattice
$$M=U\oplus \sqrt2 U\oplus (-\sqrt 2 E_8).$$

We recall that thus case is related to the
moduli space of marked Enriques surfaces, i.e. Enriques surfaces
with a choice of level 2 structure of the Picard lattice.

\smallskip

Our lattice
$$L=\sqrt 2\hbox{II}_{2,10}\cong \sqrt 2 U\oplus \sqrt2 U+\oplus (-\sqrt 2 E_8)$$
can be embedded into Kondo's lattice by means of
$$\sqrt 2 U\lo U,\quad \sqrt 2(x_1,x_2)\loma (x_1,2x_2).$$
The image consists of all integral pairs $(x_1,x_2)$ with even $x_2$.
The dual of this lattice is $(1/2)\gz\times\gz$. Hence $M'\subset L'$ corresponds
to
$$\gz\oplus\gz\oplus{1\over\sqrt2}\gz\oplus{1\over\sqrt2}\gz\oplus{-1\over\sqrt 2}E_8
\subset{1\over 2}\gz\oplus\gz
\oplus{1\over\sqrt2}\gz\oplus{1\over\sqrt2}\gz\oplus{-1\over\sqrt 2}E_8
$$
Assume that $x=\bigl(x_1,x_2,{x_3/\sqrt 2},{x_4/\sqrt 2};{\gotx/\sqrt 2}\bigr)$
is a primitive element of $M'$, which is not primitive in $L'$.
Then $x_3=2y_3$, $x_4=2y_4$ must be even and $\gotx=2\goty$ with $\goty\in E_8$.
This gives $x=(x_1,x_2,y_3\sqrt 2,y_4\sqrt2;\goty\sqrt2)$. Hence $x$ is contained
in $M$. This shows that we can apply \As.
Kondo proved that $\cz[M'/M]^{\SL(2,\gz)}$ is the direct sum of a full invariant
one-dimensional space and a 186-dimensional space $H$ which is irreducible under
$\O(M'/M)=\fz_2^8$. The elements $C\in H$ satisfy $C(0)=0$.
Hence the additive lift space of $H$ ---which is of dimension 186 by Kondo---
appears as a subspace of our 714-dimensional space:
\proclaim
{Proposition}
{The embedding of $L=\sqrt 2U\oplus \sqrt2 U\oplus (-\sqrt 2 E_8)\cong
\sqrt 2\hbox{II}_{2,10}$ into
$M=U\oplus \sqrt2 U\oplus (-\sqrt 2 E_8)$ defines an embedding of Kondo's\/
$186$-dimensional space of modular forms of weight four into our\/
$714$-dimensional space.}
KU%
\finishproclaim
\smallskip
The main result Kondo's paper is that the 186-dimensional space
defines
birational map  $\psi$  from this moduli space, i.e. $\calH_{10}/\Gamma_M$
onto its image in $P^{185}(\cz)$.
\smallskip
Unfortunately the proof is not correct, in fact
he proves that the map $\psi$ is holomorphic onto
$\calH_{10}/\Gamma_M$ and is locally biholomorphic at a special point.
But his map does not extend
to $\overline{\calH_{10}/\Gamma_M}$.
\smallskip
However Kondo's
argument works with a small modification.

In fact we know that the additive lifting from
$\cz[\fz_2^{10}]^{\SL(2,\gz)}$ to $[\Gamma_M,4]$ in injective.
Thus  we get another modular form that is relative to $\O^+(M)$,
so also in this case we add the full invariant form  to the 186
dimensional space of modular forms.\smallskip
Using the  results of proposition  \KF\ ,  it can be easily
checked that this form does not vanishes along those  0-dimensional
points that are
the base locus for the map $\psi$ . In fact these are points
parametrized by primitive isotropic elements of $M'$ that are  also in $M$.\smallskip
We
 denote by
$$A(\Gamma_M)=\sum_{r=0}^\infty[\Gamma_M,r]$$
the graded algebra of modular forms and
by
$$B(\Gamma_M)=\cz[G_0,\dots, G_{186}]\subset A(\Gamma_M)$$
the subring which is generated by the 187-dimensional
additive lift space.
The correct modification of Kondo's statement is:
\proclaim
{Theorem}
{The map
$$\overline{\calH_{10}/\Gamma_M}\to\proj(B(\Gamma_M)))\qquad(\hookrightarrow
P^{186}(\cz))$$

is everywhere regular finite and birational. The ring of modular forms
of weight divisible by 4
$$A^{(4)}(\Gamma_M):=\bigoplus_{r=0}^\infty[\Gamma_M, 4r]$$
is the normalization of the ring
$B(\Gamma_M)$ (the ring generated by the 187 dimensional additive
lift space). }
kC5
\finishproclaim

We consider now $[\Gamma_M,r]$ as subspace of $[\Gamma_L,r]$
using this concrete imbedding.
The image of $M$ in $L'/L\cong\fz_2^{12}$ is a one dimensional subspace.
We consider two six-dimensional (maximal) totally isotropic subspaces of $L'/L$ which contain
the image of $M$. The difference of their characteristic functions
defines a modular form in $[\Gamma_L,4]$. From \As\ and \PIo\ follows that this
modular form is contained in Kondo's space $[\Gamma_L,4]$.
Let now $A$ be a four dimensional totally isotropic subspace of $L'/L$ which contains the
image of $M$. From \TbK\ we obtain quadratic relations between modular forms inside
Kondo's space $[\Gamma_M,4]$.
\smallskip
We recall that in its papers Kondo's defined some quartic relations.
We don't want to go into the details of Kondo's paper and mention just that
Kondo's quartic relations are those  we
defined in \QcQ.
\proclaim
{Proposition}
{Kondo's quartic relations are a consequence of our quadratic relations.
}
RbD%
\finishproclaim
{\it Proof.\/} This follows from \QcQ.
More precisely, the model which Kondo describes in [Ko1] is nothing else
but $\proj(\calR_5/\calI_5)$.
\qed
\neupara{Relation to a result of Koike}%
There is another interesting case, which belongs to the lattice
$N=U\oplus\sqrt2U\oplus(-D_4)\oplus(-D_4)$. In this case we have $N'/N\cong\fz_2^6$.
Obviously Kondo's lattice $M$ can be embedded into $N$.
Hence we obtain enbeddings $\Gamma_L\subset\Gamma_M\subset\Gamma_N$.
The space $\cz[N'/N]^{\SL(2,\gz)}$ splits into a full invariant one-dimensional and
a 14-dimensional irreducible space under
$S_8\cong\O(\fz_2^6)$. The additive lift of the 14-dimensional space
appears as 14-dimensional subspace of Kondo's 185-dimensional space and hence
of our 714-dimensional space. We get a homomorphism
$$\calR_3/\calI_3\lo A(\Gamma_N).$$
Computer algebra shows that the dimension of $\calR_3/\calI_3$ is 6
(the dimension of the associated projective variety is 5).
The lattice $N$ admits a structure as Hermitian lattice of signature
$(1,5)$ over the ring of Gau"s integers (s.~[Ko2]).
This defines a certain 5-ball $\calB_5$
inside $\calH_{10}$. The corresponding ball quotient
is the configuration space $X(2,8)$ of 8 points in the projective line
[MY]. We obtain
a rational map
$$X(2,8)\lo \proj(\calR_3/\calI_3).$$
Koike [Koi] proved that this is a biholomorphic map.
So we recover Koike's observation that certain quartic relations are
consequences of quadratic ones and even more that all these relation
live already in the ten-dimensional space $\calH_{10}/\Gamma_L$.

\neupara{Final Remark}%
in \RFb\ we used the
full 715-dimensional additive lift space. But the relations which we described
play in the 714-dimensional irreducible part.
It may be true that the invariant form is not necessary,
because it may be possible that the second symmetric
power of the 714- and 715-dimensional
space is the same. We don't know whether this is true or not.
\smallskip
The invariant form in the
715-dimensional space possible being superfluously does not
detain it to be in some
sense very basically. To explain this we consider a rational orthogonal
transformation $g\in\O^+(\hbox{II}_{2,10}\otimes_\gz\qz)$ with the property
$$g\O^+(\hbox{II}_{2,10}[2])g^{-1}\subset \O^+(\hbox{II}_{2,10}).$$
Then one can transform $f$ with $g$ to obtain a form on $\O^+(\hbox{II}_{2,10})[2]$.
To give a simple example we realize $\hbox{II}_{2,10}$ as
$$\hbox{II}_{2,10}=\hbox{II}_{1,1}\times\hbox{II}_{1,9}.$$
We take for $g$ the orthogonal transformation
$$(a,b,C)\loma (a/2,2b,C).$$
To understand it better we consider the standard embedding of $\calH_{10}$
into \quad $\hbox{II}_{1,9}\otimes_\gz\cz$:
The modular form $f$ is determined by the function
$$F(Z):=f(1,*,Z),$$
where the star is taken such that the entry is in the zero-quadric.
Now the effect on $F$ is given by the transformation $F(Z)\mapsto F(2Z)$.
It is possible to compute the values at the cusps and to prove in this way that
the transformed form is contained in the 715-dimensional additive lift space.
This  gives:
\proclaim
{Proposition}
{The $715$-dimensional additive lift space (\Sev) in the space of
modular forms $[\O^+(\hbox{II}_{2,10})[2],4]$
is generated as $\O(\fz_2^{10})$-module by the single form $F(2Z)$, where
$F(Z)$ denotes the full invariant form, expressed in standard coordinates.}
Vea%
\finishproclaim
\need6cm
\vskip2cm\noindent
{\paragratit References}\medni
\medskip
\item{[AF]} Allcock,\ D.\ and Freitag,\ E.:
{\it Cubic surfaces and Borcherds products},
Commentarii Math.\ Helv. Vol.~{\bf 77}, Issue 2, 270--296 (2002)
\medskip
\item{[BB]} Baily,\ W.L., and Borel, A.:
{\it Compactification of arithmetic quotients of bounded
Symmetric domains,}  Annals of Math. {\bf 84},
No 3, 442--528  (1966)
\medskip
\item{[Bo1]} Borcherds,\ R.: {\it Automorphic forms
with singularities on Grassmannians,}
Invent. math. {\bf 132}, 491--562 (1998)
\medskip
\item{[Bo2]} Borcherds,\ R.: {\it The Gross-Kohnen-Zagier
theorem in higher dimensions, }
Duke Math.~J. {\bf 97}, No 219--233 (1999)
\medskip
\item{[BK]} Brunier, J. Kuss, M.: {\it Eisensetin series attached to lattices and modular
forms on orthogonal groups,\/} Manuscr.\ Math.~{\bf 106}, 443--450 (2001)
\medskip
\item{[CS]} J.\ H.\ Conway, N.\ J.\ A.\ Sloane.
Sphere packings, lattices
and groups. Springer-Verlag New York,
Grundlehren der mathematischen
Wis\-sen\-schaf\-ten {\bf 290}, (1988)
\medskip
 \item{[Fa]} Fay, J.: {\it On the Riemann-Jacobi formula,\/}
 Nachr. Akad. Wiss. Gštt., II. Math.-Phys. Kl. 1979 61-73 (1979)
\medskip
\item{[FH]} Freitag, E. Hermann, C.F.: {\it Some modular
varieties in low dimension,\/} Advances in Math. {\bf 152}, 203-287 (2000)
\medskip
\item{[Fr1]} Freitag, E.: {\it Modulformen zweiten Grades zum rationalen und Gau"sschen
Zahlk"orper,\/} Sitzungsberichte der Heidelberger
Akademie der Wis\-sen\-schaf\-ten, 1.~Abh.
(1967)
\medskip
\item{[Fr2]} Freitag, E.: {\it Some modular
forms related to cubic surfaces\/}
Kyungpook Math. J. {\bf 43}, No.3, 433-462 (2003)
\medskip
\item{[Fr3]} Freitag, E.:
{\it Comparison of different models of the moduli space
of marked cubic surfaces,\/}
Proceedings of Japanese-German Seminar, Ryushi-do, edited by T.~Ibukyama
and W.~Kohnen, 74-79 (2002)
\medskip
\item{[FS]} Freitag, E. Salvati-Manni, R.:
{\it Some modular varieties of low dimension II,\/} preprint (2004)
\medskip
\item{[Ig]} Igusa, J.:
{\it On the graded ring of theta-constants \/}
Am. J. Math. {\bf 86},, 219-246 (1964).
\medskip
\item{[Kn]} Kneser, M.: {\it Erzeugung ganzzahliger orthogonaler Gruppen
durch Spiegelungen,\/} Math.~Ann.~{\bf 255}, 453--462 (1981)
\medskip
\item{[Ko1]} Kondo, S.: {\it The moduli space of
Enriques surfaces and Borcherds products,\/}
J. Algebraic Geometry  {\bf 11}, 601-627(2002)
\medskip
\item{[Ko2]} Kondo, S.: {\it  The moduli
space of 8 points on $P^{1}(\cz)$ and automorphic forms\/}, preprint (2005)
\medskip
\item{[Koi]} Koike, K.: {\it  The projective
embedding of the configuration space $X(2,8)$\/}, preprint (2005)
\medskip
\item{[MY]} Matsumoto, K. Yoshida, M.: {\it Configuration space of 8
points on the projective line and a 5-dimensional Picard modular group,\/}
Compositio Mtah.~{\bf 86}, 265--280 (1993)
\medskip
\item{[Ni]} Nikulin,\ V.V.:
{\it Integral symmetric bilinear forms and
some of their applications,\/}
Math.\ USSR Izvestija {\bf 14}, No 1    (1980)
\medskip
\item{[Sch]} Scheithauer, N.: {\it Moonshine for Conway's group,\/}
Habilitationsschrift, Heidelberg (2004)
\bye